\documentclass[11pt]{article}
\usepackage{mathtools}
\usepackage{amsmath,amsthm,amssymb}
\usepackage{hyperref}
\usepackage{footnote}
\usepackage{tablefootnote}
\usepackage{lipsum}
\usepackage{siunitx}
\usepackage{authblk}

\begin{document}
%
%
%

\author{Srinibas Swain\footnote{ssrinibas@gmail.com, Indian Institute of Information Technology Guwahati, Assam 781015, India.}, C. Paul Bonnington\footnote{paul.bonnington@monash.edu, Monash eResearch Centre, Monash University, Clayton, Victoria 3800, Australia.}, Graham Farr\footnote{graham.farr@monash.edu, Faculty of Information Technology, Monash University, Clayton, Victoria 3800, Australia.}, Kerri Morgan\footnote{kerri.morgan@deakin.edu.au, Deakin University, Geelong, Australia, School of Information Technology, Faculty of Science Engineering \& Built Environment, Australia.}}

\title{A survey of repositories in graph theory}
\maketitle
\begin{abstract}
Since the pioneering work of R. M. Foster in the 1930s, many graph repositories have been created to support research in graph theory. This survey reviews many of these graph repositories and summarises the scope and contents of each repository. We identify opportunities for the development of repositories that can be queried in more flexible ways.
\end{abstract}
\section{Introduction}
The history of graph theory may be traced to 1736, when the Swiss mathematician Leonhard Euler solved the K\"onigsberg bridge problem~\cite{RD2010}. In recent decades, graph theory has established  itself as an important mathematical tool in a wide variety  of subjects, including computer science, electrical engineering, operational research, chemistry, genetics, linguistics, geography, sociology and architecture~\cite{TWG}. At the same time it has also become established as a sophisticated mathematical  discipline in its own right. 

Readily available graph data may help graph theory researchers to develop insight and intuition. However generating and storing graph data is extremely challenging due to its large magnitude. The number of graphs grows exponentially with increasing number of vertices~\cite{HP1973,RW1998}. The numbers of labelled and unlabelled connected graphs of order $<13$ are listed in Table~\ref{t1}.  
\begin{table}[htp]
\begin{tabular}{ |r|r|r| } 
\hline
\# vertices & \# connected labelled graphs & \# connected unlabelled graphs \\
\hline
$1$ & $1$ & $1$\\
\hline
$2$ & $1$ & $1$\\
\hline
$3$ & $4$ & $2$\\
\hline
$4$ & $38$ & $6$\\
\hline
$5$ & $728$ & $21$\\
\hline
$6$ & $\num{26704}$ & $112$\\
\hline
$7$ & $\num{1866256}$ & $853$\\
\hline
$8$ & $\num{251548592}$ & $\num{11117}$\\
\hline
$9$ & $\num{66296291072}$ & $\num{261080}$\\
\hline
$10$ & $\num{34496488594816}$ & $\num{11716571}$\\
\hline
$11$ & $\num{35641657548953344}$ & $\num{1006700565}$\\
\hline
$12$ & $\num{73354596206766622208}$ & $\num{164059830476}$\\
\hline
\end{tabular}
\caption{Growth of number of graphs~\cite{HP1973, RW1998}.}\label{t1}
\end{table}

A \emph{graph parameter} is defined as a function $f:\{\hbox{graphs}\}\rightarrow \mathbb{Z}$ that is invariant under isomorphism. We sometimes extend the term to include non-integer-values like complex numbers and polynomials. Some examples of graph parameters are chromatic index, treewidth, and the Tutte polynomial. Finding the chromatic index is known to be NP-complete~\cite{GJ1979}. The chromatic index is useful in scheduling problems~\cite{JAN} and routing problems~\cite{JA}. Treewidth was introduced by Robertson and Seymour in 1983~\cite{RS}. Treewidth is discussed in detail in~\cite{BD}. Treewidth plays an important role in parameterized complexity~\cite{BD}. Many graph-theoretic problems which are NP-complete are polynomial time solvable when restricted to graphs of bounded treewidth~\cite{BD}. The problem of finding the treewidth of a graph is NP-complete~\cite{GJ1979}. The Tutte polynomial is a two-variable polynomial introduced by W. T. Tutte~\cite{wttu}. It plays a key role in the study of counting problems on graphs, and has close connections with statistical mechanics and knot theory~\cite{DW, Welsh}. 

A repository of graphs and their related parameters provides a framework that enables researchers to test new theorems and conjectures and may also facilitate new mathematical discoveries related to the data present in the repository.

Some graph parameters are computable in polynomial time, some take exponential time and some of them are uncomputable~\cite{AR2009,GJ1979}. There are thousands of graph parameters that are of interest to researchers. This makes the existence of any substantial graph repository that contains all published graph  parameters impractical. Several researchers have created useful graph repositories for specific graph parameters. Some of the major challenges in creating these repositories are computing the parameters, representing them and storing the data. Therefore, these repositories can be broadly classified as repositories before and after the advent of magnetic media. 

We list the graph repositories in the order of their creation in Section~\ref{GR}. The graph repositories that use the printed form to store the graphs are discussed in Section~\ref{prr}. The repositories which store the graphs in electronic media are listed in Section~\ref{err}. The interactive graph repositories are listed in Section~\ref{irr}. We summarise the sources of graphs for the repositories and their interdependencies in Section~\ref{nsh}. We conclude the survey by outlining the key features for some of these repositories (listed in Section~\ref{GR}) in Section~\ref{conc}. 

Note that the graph repositories discussed in this survey are specific to research in graph theory and are different to the graph repositories used for running experiments and tests on graph algorithms, such as social networks, biological networks, collaboration networks and ecological networks~\cite{nr}.

\section{Graph repositories}
\label{GR}
In this section, we discuss the various commonly used graph repositories. We discuss the parameters that each of these repositories contains and some of their limitations. Note that we associate the graph information described in this section with the time at which it appeared in its repository, although in some cases the process of generating these data might have started prior to the time the data appeared in the repositories.

Table~\ref{t1} shows that storing all graphs with order greater than $12$ is challenging. To deal with this storage limitation, several researchers prefer to generate graphs on an as-needed basis. Therefore, graph generators can also regarded as a rich source of graph repositories. However, in this survey we focus more on repositories that store graphs, after briefly discussing graph generators in the next subsection.
\subsection{Graph generators}\label{gg}
In this section we highlight various graph generators. Some classes of graphs are complex to generate, for example strongly-regular graphs, snarks, non-Hamiltonian graphs, and colour-critical graphs, whereas simple graphs, bipartite graphs, regular graphs and trees are relatively easy to generate. For the former class of graphs it is important to store the graphs to make the most of the computational effort of generating them, but for the latter efficient generators can be more useful.   
\subsubsection{Initial attempts (1946--1976)}
In this section we list early work on graph generation. The first such attempt was reported by Kagno~\cite{kagno} in 1946. 
In this work the focus was to determine whether a given graph of degree $\le 6$ has a non-trivial automorphism group. He generated the $156$ undirected graphs of order $6$.
Subsequently Heap~\cite{heap}, Baker et al.\ \cite{dewdney} generated the $12346$ and $274668$ graphs with order $8$ and $9$ respectively. Similarly Read~\cite{rread} generated the $9608$ digraphs of order $5$.\ Morris~\cite{morris} and Frazer~\cite{frazer} generated the $1301$ and  $123867$ trees of order $13$ and $18$ respectively.\ McWha~\cite{mcwha} generated the $456$ tournaments of order $7$. Morris~\cite{morriss} generated the $10$ and $36$ self-complementary graphs of order $8$ and $9$ respectively. 
Bussemaker et al.~\cite{busmaker} generated the $509$ connected cubic graphs of order $14$.

Understandably, all the works mentioned above did not report generation of graphs of 
orders higher than the aforementioned values, due to limitations on computation~\cite{RR81}. However, these attempts were useful as they were the first ones to use computers to generate graphs and later provided information on designing algorithms for efficient graph generation. There are many attempts at graph generation made from 1976 until now. In the following two subsections we briefly describe two of those which are versatile in nature and commonly used. Other generators are mentioned in the Section~\ref{err}. 
\subsubsection{Chemical \& abstract Graph environment (CaGe) (1997--present)}
CaGe was developed by G. Brinkmann, O. D. Friedrichs, S. Lisken, A. Peeters, N. V. Cleemput~\cite{cage}. CaGe can be accessed using the link \url{ http://caagt.ugent.be/}. It is used to generate graphs, often those relate to interesting chemical molecules. CaGe consists of generators for following classes of graphs: $k$-regular plane ($k\le 4$), hydrocarbons, tubes and cones, triangulations, quadrangulations, and general plane. 
CaGe provides the user to view graphs in different formats: three-dimensional, two-dimensional and adjacency information. 
\subsubsection{SageMath (2005--present)}
SageMath (previously Sage or SAGE) is an interface for studying pure and applied mathematics~\cite{sagemath}. This covers a variety of topics in mathematics including algebra, calculus, number theory, cryptography, combinatorics, and graph theory. In this section we discuss the graph theory aspect.

SageMath provides a framework to generate simple graphs, digraphs, random graphs, hypergraphs. 
Sage Graphs can be created from a wide range of input formats, for example, graph6, sparse6, adjacency matrix, incidence matrix, and list of edges.

It also incorporates a database of strongly regular graphs, and it returns a specific strongly regular graph based on a user's request, when one exists. This uses Andries Brouwer's database (\url{https://www.win.tue.nl/~aeb/graphs/srg/srgtab.html}) of strongly regular graphs to return non-existence results~\cite{cp17}. 

Although SageMath is a graph generator, its developers maintain a database (repository) of all unlabeled graphs with order at most $7$. SageMath provides the facility of querying for graphs (using the aforementioned database) that satisfy constraints on a certain set of parameters, including numbers of vertices and edges, density, maximum and minimum degree, diameter, radius, and connectivity.

SageMath also contains many algorithms to compute different graph parameters, for example chromatic number, genus, Tutte polynomial and to check if a graph is asteroidal triple-free\footnote{An independent triple $\{x, y, z\}$ is called an asteroidal triple (AT, for short) if between any pair in the triple there exists a path that avoids the neighbourhood of the third vertex.}. The full list of the algorithms is available at \url{http://doc.sagemath.org/html/en/reference/graphs/index.html}. 

The generators in Section~\ref{gg} provide users with the facility to generate data using 
their local computing facility. These are useful if the user's aim is only to generate graphs supported by these generators. However these generators do not facilitate searching for graphs satisfying various graph properties. 

\subsection{Print resources}\label{prr}
The following repositories of the late twentieth century used printed form to store information. A large collection of graph data is very difficult to present in this form because of the inherent limitations of print media. Despite these limitations, these repositories were quite significant for specific problems and very well structured. Although the repositories listed in this section are quite rich in content, it is hard to search, insert, update or extend data in these type of repositories. These repositories do not have the capability to handle queries which look up some range of parameter values or combination of several parameters, e.g., list all connected graphs whose diameter is between $4$ and $7$ with chromatic index $5$.
\subsubsection{Data compendium of linear graphs (1967)}\label{dcl}
\begin{figure}
 \includegraphics[width=\linewidth]{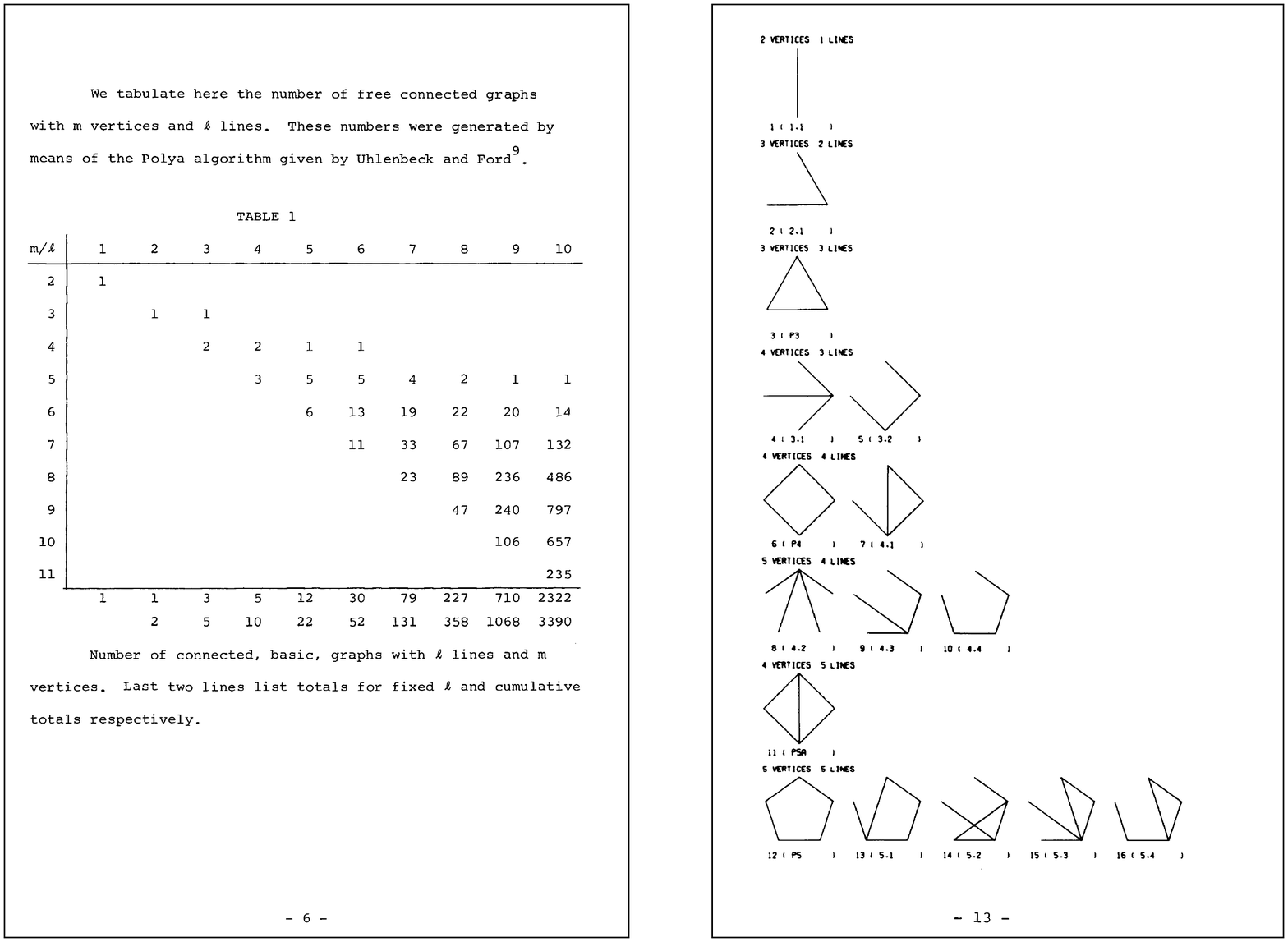} 
\caption{An extract from the 1967 report by Baker et al.~\cite {bak1967}.}
 \label{bak}
\end{figure}

In 1967, Baker, Gilbert, Eve, and Rushbrooke published a list of $1460$ simple connected sparse graphs with number of edges $\leq 10$ and order $\leq 11$ ~\cite{bak1967}. They generated these graphs for a project aimed at computing high-temperature expansions for the field dependent free energy of a Heisenberg model ferromagnet. It took three years of computation to accumulate this data. Although this list of graphs reported by Baker et al.\ is not the complete set of graphs of order $\leq 11$, it still had lot of practical significance in the field of cooperative phenomena~\cite{CD1960,PH1964}. An example of the data presented in \cite{bak1967} is given in Figure~\ref{bak}. 


Given the resources available at that time this repository, although small, was a valuable resource and one of the first few repositories of graphs. However, this repository does not store any graph parameters.
\subsubsection{The Foster Census (1930--1988)}\label{tfc}
\begin{figure}
\center
 \includegraphics[width=.8\linewidth]{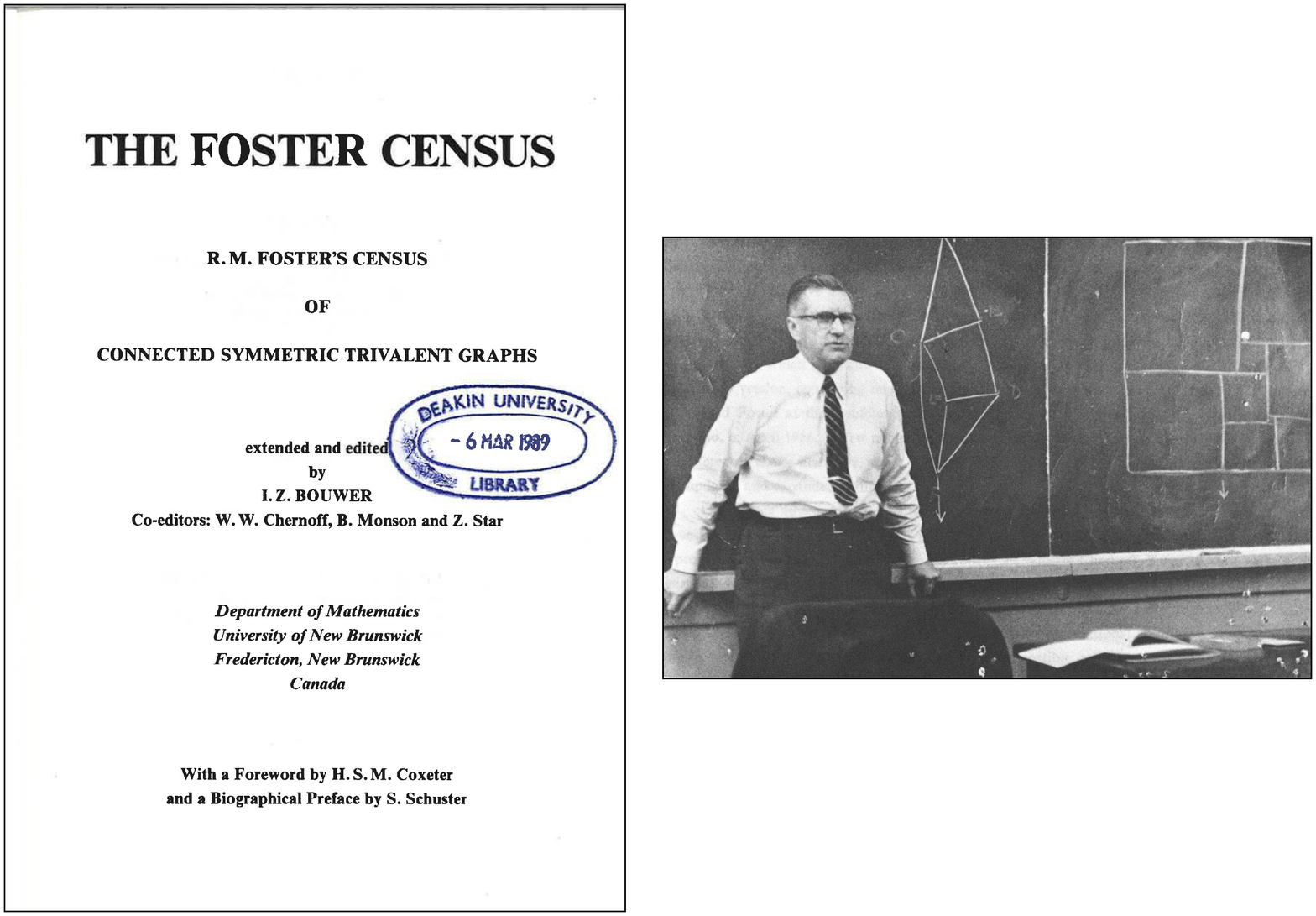} 
\caption{Foster's Census with Foster, from~\cite{FB1988}.}\label{Foster}
\end{figure}

Although the data compendium by Baker et al.\ was published in 1967, Ronald M. Foster started storing graph data decades before  the compendium was published.\ In the 1930s, Foster started collecting small cubic symmetric graphs while he was employed by Bell Labs~\cite{FB1988}. Foster's hand-prepared list was surprisingly accurate and had only one omission for graphs of order up to 240. On the other hand, for order between 240 and 512 there
were several other omissions. In $1988$, when Foster was 92, the Foster Census listing all cubic symmetric graphs up to 512 vertices was compiled by I. Z. Bouwer, W. W. Chernoff, B. Monson and Z. Star (``The Foster Census"~\cite{FB1988}). The cover page of the Foster Census and a picture of Foster is shown in Figure~\ref{Foster}.\\

\noindent{\bf The missing graphs and the discovery}\\
Conder and Dobcs\'anyi~\cite{MC2000} found all the missing graphs from the Foster Census~\cite{FB1988} while they were compiling all cubic symmetric graphs with order up to $768$.

There was a graph among these missing graphs (448C, as per Foster’s convention)
which was the smallest graph that is arc-transitive but has no involutory automorphism
reversing an arc. Interestingly, the previous smallest known graph of this type was a graph
with 6652800 vertices~\cite{MC2000}. This illustrates the research benefits of compiling and using graph repositories. 


The first version of the Foster census not only had the list of cubic symmetric graphs, but also stored some useful graph parameters like diameter, girth, $s$-transitivity, hamiltonicity and bipartiteness of each listed graph. 
\subsubsection{An Atlas of Graphs (1998)}\label{aag}
As we noted in the previous sections, all published graph repositories before the advent of magnetic media were narrow in scope. The first comprehensive graph data resource was Read and Wilson's ``An Atlas of Graphs"~\cite{RW1998}. This repository has all graphs with up to $7$ vertices and all digraphs with up to $5$ vertices. For each graph, Read and Wilson included values of some parameters and also gave information on whether or not certain properties hold. The Atlas also contains some families of graphs. We list some of the graphs and graph parameters from the book.
\begin{itemize}
\item {\bf Simple graphs and parameters}\\ This book not only lists all the unlabelled simple graphs of order $\leq 7$ but also depicts the graphs. The graphs are listed in increasing order of number of vertices. For a fixed number of vertices the graphs are listed in increasing order of edges. For fixed numbers of vertices and edges, the graphs are listed in increasing order of degree sequence.

For each of the aforementioned graphs the book lists the following graph parameters: order, number of edges, degree sequence, number of connected components, girth, number of cycles of shortest length, diameter, clique number, independence number, vertex connectivity, edge connectivity, number of automorphisms, complement, whether or not certain properties hold (like bipartiteness, Eulerian, forest, Hamiltonian, planar, tree, uniquely colourable), chromatic number, chromatic index, chromatic polynomial and spectral polynomial.

\item {\bf Trees} \\The book lists trees, rooted trees, homeomorphically irreducible trees\footnote{A tree in which all nodes have degree other than 2 is called a {\it homeomorphically irreducible tree}.} and identity trees of order $\leq 30$. It also depicts the $987$ trees with up to $12$ vertices (together with their degree sequences), homeomorphically irreducible trees with up to $16$ vertices, the identity trees up to $14$ vertices and binary trees with up to $7$ vertices. It also lists some of the parameters of trees: order, number of edges, degree sequence, diameter, independence number, number of automorphisms of the graph, properties (like bicentral, bicentroidal, central, centroidal, homeomorphically irreducible, identity tree) and spectral polynomial.

\item {\bf Regular graphs} \\It lists the number of labelled cubic and connected cubic graphs up to $40$ vertices and 4-regular graphs up to $15$ vertices. It also depicts the connected cubic graphs  with up to $14$ vertices, 4-regular graphs up to $11$ vertices, and 5-regular and 6-regular graphs up to $10$ vertices. It also depicts the connected bicubic graphs with up to $16$ vertices and the cubic polyhedral graphs (without triangles) with up to $18$ vertices, connected cubic transitive graphs with up to $34$ vertices, 4-regular transitive graphs with up to $19$ vertices and symmetric graphs with up to $54$ vertices. The book also lists all the parameters for the regular graphs which are listed for simple graphs except the graph polynomials.

\item{\bf Other classes of graphs}\\
The scope of this book with respect to variety of graphs is wide. This book lists the number of bipartite graphs, connected bipartite graphs, unicyclic graphs and self-complementary\footnote{A self-complementary graph is one that is isomorphic to its complement.} graphs with up to $20$ vertices. It also listed even graphs\footnote{ A connected graph $G$ is called {\it even} if for each vertex $v$ of $G$ there is a unique vertex $u$ such that $d(v,u) = \hbox{diam } G$}, Eulerian graphs, and connected line graphs up to $16$ vertices. This lists all Hamiltonian graphs up to $11$ vertices. 

The book also depicts all connected bipartite graphs with up to $8$ vertices, Eulerian graphs with up to $8$ vertices, self-complementary graphs up to $9$ vertices, connected triangle-free graphs up to $10$ vertices (none of degrees less than $3$), connected line graphs up to $8$ vertices and unicyclic graphs up to $8$ vertices.

\item{\bf Planar graphs}\\
The authors depict planar graphs as proper planar embeddings marking the interior and exterior faces. The book depicts 2-connected plane graphs up to $7$ vertices and $3$-connected graphs with up to $8$ vertices together with their respective degree sequences. The book also depicts the outerplanar graphs with up to $9$ vertices.

\item{\bf Digraphs}\\
The book is a repository of digraphs as well. It lists various types of digraphs (connected, unilateral, strong, acyclic, self-complementary, self-converse and tournaments) with up to 11 vertices. The digraphs are listed in increasing order of number of vertices. For a fixed number of vertices, they are listed in increasing order of number of arcs. For fixed numbers of vertices and arcs, digraphs are listed in increasing order of degree sequences. For fixed numbers of vertices and arcs, and fixed degree sequences, they are listed in order of increasing number of automorphisms.

The book depicts the acyclic digraphs with up to $5$ vertices. It also depicts Eulerian digraphs up to $5$ vertices, $2$-regular digraphs up to $7$ vertices, self-complementary digraphs up to $5$ vertices and tournaments up to $7$ vertices. It also differentiates the strong tournaments among the depicted tournaments. Weakly connected transitive digraphs up to $5$ vertices are also depicted. 

Along with the digraphs the book also lists some important graph parameters associated with them. The parameters listed for digraphs are: number of vertices, number of arcs, out-degree sequence and in-degree sequence, connectivity (whether digraph is disconnected, connected but not unilateral, unilateral but not strong or strong), number of automorphisms and whether or not certain properties hold (acyclic, Eulerian, Hamiltonian, self-complementary, tournament, self-converse).

\item{\bf Signed graphs}\\
This book is one of the few repositories that lists and depicts signed graphs. It lists signed graphs, balanced signed graphs and signed trees with up to $12$ vertices. The book depicts signed graphs with up to $5$ vertices and signed trees with up to $7$ vertices. 

The book lists parameters including number of vertices, number of edges (numbers of positive and negative edges), number of $3$-cycles, number of $4$-cycles, and number of $5$-cycles, and also specifies whether the signed graph is balanced or unbalanced for all the signed graphs.
\item{\bf Graphs and Ramsey numbers}\\
Ramsey numbers play an important role in graph theory and in probabilistic methods. This book gives the Ramsey numbers for some pairs of connected graphs with up to $5$ vertices. It also depicts all isolate-free graphs\footnote{ A graph
is {\it isolate-free} if it has no isolated vertex.} with up to $7$ edges, with their Ramsey numbers. Each depiction also contains the serial number of the graph in Burr's catalogue~\cite{burr}. The authors also depict graphs with more than $7$ edges for which the Ramsey number is known.

\item{\bf Polynomials}\\
Graph polynomials can be useful as some of them encapsulate a lot of information about the graph. Read and Wilson listed the chromatic polynomials of graphs (up to $7$ vertices), cubic graphs (up to $14$ vertices), and 4-regular graphs (up to $11$ vertices). Each polynomial is presented in the \emph{power form} (in decreasing powers of $\lambda $), in sum of falling factorial form, and in \emph{tree form}, which for a graph with $k$ components, is the polynomial $P$ such that $\lambda^k P(\lambda - 1)$ is the chromatic polynomial.\footnote{Read and Wilson~\cite{RW1998} give the example of the chromatic polynomial of the butterfly graph $G$, with $V(G)=\{0,1,2,3,4\}$ and $E(G)=\{01,02,03,04,12,34\}$, presented in the different forms:\begin{itemize}
\item Power form: $\lambda^{5}-6 \lambda^{4} +13 \lambda^3-12\lambda^2+4\lambda$
\item Falling factorial form: $\lambda^{(5)}+\lambda^{(4)}+\lambda^{(3)}$, where $\lambda^{(r)}=\lambda(\lambda-1)\cdots(\lambda-r+1)$ 
\item Tree form: $\lambda P(\lambda - 1)=\lambda((\lambda-1)^4-2(\lambda-1)^3+(\lambda-1)^2)=\lambda(\lambda-1)^4-2\lambda(\lambda-1)^3+\lambda(\lambda-1)^2$,
so $P(x)=x^4-2x^3+x^2.$
\end{itemize}} 

The authors also list the characteristic polynomials of graphs (up to $7$ vertices), trees (up to $12$ vertices), cubic graphs (up to $14$ vertices) and 4-regular graphs (up to $11$ vertices). For graphs and trees, each polynomial is presented in power form, together with its frequency and spectrum (usually with zero and integer eigenvalues listed first, and the remainder listed in decreasing order); for cubic and 4-regular graphs only the polynomial is given.
\item{\bf Special graphs}\\
There are some graphs which are considered special due to some specific properties that these graphs satisfy. They play an important role in theorem and conjecture verification. The authors of this book lists some of them. The special graphs listed and depicted in this book are platonic graphs, Archimedean graph, M\"obius graph, cages, non-Hamiltonian cubic graphs, generalised Petersen graphs, snarks, graphs drawn with minimum crossings, the two smallest cubic identity graphs, the hypercube, the Greenwood-Gleason graph, cubic graphs with no perfect matching, the Goldner-Harary graph, the Biggs-Smith graph, Folkman's graph, Tietze's graph, Meredith's graph, Chv\'atal's graph, Franklin's graph, the Moser spindle, the Herschel graph, Mycielski's graph or Gr\"otzsch's graph, and Royle's graph. All these graphs are defined in Chapter~6 of the book~\cite{RW1998}.
\end{itemize} 
\begin{table}[htp]
\centering
\begin{tabular}{ |l|c|c| } 
\hline
Graph type & Order & Exhaustive list (Y/N)\\
\hline
\hline
Unlabelled simple& 1--7 & Y \\
\hline
Trees  & 1--12 & Y \\
\hline
Homeomorphically irreducible trees  & 1--16 & Y \\
\hline
Identity trees  & 7--14 & Y \\
\hline
Binary trees  & 1--7 & Y \\
\hline
Connected cubic& 4--14 & Y \\
\hline
Connected 4-regular& 5--11 & Y \\
\hline
Connected 5-regular  &  6--10 & Y \\
\hline
Connected 6-regular  & 7--10&Y  \\
\hline
Connected bicubic  & 4--16 & Y \\
\hline
Connected polyhedral  & 8--18 &  Y\\
\hline
Connected cubic transitive&4--34& Y \\
\hline
Connected 4-regular transitive& 5--19 & Y \\
\hline
Symmetric cubic  & 4--54 & Y \\
\hline
Connected bipartite & 2--8 & Y\\
\hline
Eulerian & 1--8 & Y\\
\hline
Self-complementary & 4--9 & Y\\
\hline
Connected triangle-free & 6--10 & Y \\ \hline
Connected line & 1--8 & Y\\ \hline
Unicyclic & 3--9 & Y\\ \hline
Plane ($2$-connected, $3$-connected) & 3--7, 4--8 & Y \\ \hline
Outerplanar & 3--9  & Y \\ \hline
Digraphs & 1--4 & Y \\ \hline
Acyclic digraphs & 1--5  & Y \\ \hline
Eulerian digraphs & 1--5 & Y \\ \hline
2-regular digraphs & 3--7 & Y \\ \hline
Self-complementary digraphs& 1--5 & Y \\ \hline
Tournaments & 1--7 & Y \\ \hline
Weakly connected transitive digraphs&1--4 & Y \\ \hline
Signed& 1--5 & Y \\ \hline
Signed trees& 1--7 & Y \\ \hline
Other known graphs & -- & N \\ \hline
\end{tabular}
\caption{A summary of graphs in ``An Atlas of Graphs".}\label{tga}
\end{table}
\begin{table}[h]
\centering
\begin{tabular}{ |l|c| } 
\hline
Type of graph(s) & List of parameters\\ \hline
\hline
Unlabelled graph, tree, & automorphism group size, diameter, independence number,\\
Regular graph & degree sequence, spectral polynomial \\ \hline

Unlabelled graph, & girth, circumference, tree, clique number, chromatic number,\\
Regular graph & bipartiteness, chromatic index, chromatic polynomial,\\
 & uniquely colourable, Eulerianness, planarity, Hamiltonicity,\\
 & number of shortest-length cycles, \\
 &   vertex-connectivity, edge-connectivity\\ \hline

Unlabelled graph & number of components, tree \\ \hline

Tree & bicentral, bicentroidal, central, centroidal, identity tree,\\ & homeomorphically irreducible\\ \hline

Digraph & acyclic, disconnected, connected but not unilateral,\\
&  unilateral but not strong, strong, Eulerianness, Hamiltonicity, \\
& degree sequence (out-degree and in-degree), automorphism group \\
& size, self-complementary, tournament, self-converse \\ \hline

Signed graph & number of edges (positive, negative), 3-cycle,\\
& 4-cycle, 5-cycle, balanced or unbalanced \\ \hline
\end{tabular}
\caption{A summary of parameters in ``An Atlas of Graphs"}. \label{pga}
\end{table}
We discussed the limitations of print media in Section~\ref{prr}. ``An Atlas of Graphs" is the most comprehensive printed graph repository. A summary of the data in ``An Atlas of Graphs" is given in Tables~\ref{tga} and \ref{pga}. These tables demonstrate the richness of this repository. Searching for a graph that the book contains is straightforward. However, constructing a collection of graphs satisfying some specific conditions can be a laborious task with any printed repository. If the data present in ``An Atlas of Graphs" is represented in electronic media, searching and constructing data would be easier and efficient.
\subsection{Electronic resources}\label{err}
The limitations of printed form, and the challenges in generating and storing data for higher order graphs, prompted researchers to build electronic representations of graph data that can be accessed easily. 

As discussed, one of the major challenges in creating a graph repository is the representation of graphs. A straightforward way to store a graph in electronic form is via its adjacency matrix, which requires $n^2$ bits to store a graph of order $n$.\ McKay introduced {\tt graph6} and {\tt sparse6} formats for storing undirected graphs in a compact manner~\cite{MK1984}, using only printable ASCII characters obtained from bit manipulation of the adjacency matrix. Files in these formats have text type and contain one line per graph. The format {\tt graph6} is suitable for small graphs, or large dense graphs, while the format {\tt sparse6} is more space-efficient for large sparse graphs.     

Many graph theorists generate various sets of graph data in the course of their research. Substantial contributions of this type been made by Conder~\cite{MC2002}, McKay~\cite{MK19} and Royle~\cite{GR2007}.
\subsubsection{Brendan McKay's combinatorial data (1984--present)}\label{bdm}
McKay maintains a collection of graphs, latin squares, cubes and Hadamard matrices. In this section we will only discuss the graph-related data. McKay's collection is one of the largest available repositories of graphs. It contains more than 150 million graphs in total~\cite{MK19}. McKay's collection of graphs is interesting and useful for its variety and magnitude. We list various graphs from McKay's repository below. 

\begin{itemize}
\item {\bf Simple graphs}\\
McKay stored all connected unlabelled graphs of order $\leq 11$, and generated all such graphs of order $\leq 13$.  McKay generated these graphs according to their numbers of edges and vertices. The program used to generate these graphs is {\tt geng} \cite {MK19}.
\item {\bf Special class of graphs}\\
McKay's repository contains all Eulerian graphs (and all connected Eulerian graphs) and chordal graphs up to 12 vertices. 

McKay's repository contains all perfect graphs up to 11 vertices.

Let $G$ be a regular graph with $n$ vertices and degree $k$. $G$ is said to be {\it strongly regular} if there exist integers $\lambda$ and $\mu$ such that:
\begin{itemize}
  \item every two adjacent vertices have $\lambda$ common neighbours, and
  \item every two non-adjacent vertices have $\mu$ common neighbours.
    
\end{itemize}
A graph of this kind is sometimes said to be an SRG$(n,k,\lambda, \mu)$. Strongly regular graphs were introduced by Raj Chandra Bose in 1963~\cite{RAJ63}. McKay's collection of strongly regular graphs is one of the most comprehensive lists of strongly regular graphs. Most of these graphs have been computed by McKay and/or Ted Spence. One of the most significant results by McKay and Spence is the classification of regular two-graphs on $36$ and $38$ vertices~\cite{MS}. An immediate consequence of this was that all strongly regular graphs with parameters (35, 16, 6, 8), (36, 14, 4, 6), (36, 20, 10, 12) and their complements are known. 

A graph is \emph{hypohamiltonian} if it is not Hamiltonian but each graph obtained from it by removing one vertex is Hamiltonian. Petersen graph is the smallest (order) hypohamiltonian graph. Table~\ref{mga} summarises McKay's collection of hypohamiltonian graphs.

All non-isomorphic connected planar graphs with up to 11 vertices are stored in McKay's collection of graphs. McKay also stores the planar embeddings of the graphs.
McKay's repository stores plane 5-regular simple connected graphs up to 36 vertices and nonhamiltonian planar cubic graphs (this has graphs with no faces of size 3,  cyclically 4-connected graphs, graphs with no faces of size 3 or 4 with cyclic connectivity exactly 4, and cyclically 5-connected graphs) and hypohamiltonian planar graphs (this includes cubic graphs of girth 4, cubic planar graphs of girth 5, and cubic planar graphs with an \emph{$\alpha$-edge}, where an $\alpha$-edge in a graph is an edge which is present on every Hamiltonian cycle).

Self-complementary graphs can have only orders congruent to 0 or 1 modulo 4. McKay stores all such graphs up to 17 vertices. However, he has a partial list of graphs for 20 vertices. McKay stores $\num{8571844}$ self-complementary graphs (for 20 vertices) out of  $\num{9168331776}$ graphs.

A connected graph is \emph{highly irregular} if the neighbours of each vertex have distinct degrees. Such graphs exist for all orders except 3, 5 and 7. All highly irregular graphs with up to 15 vertices are listed in McKay's combinatorial data.

\item{\bf Ramsey graphs}:\\
A \emph{Ramsey$(s,t,n)$-graph} is a graph with $n$ vertices, no clique of size $s$, and no independent set of size $t$. A \emph{Ramsey$(s,t)$-graph} is a Ramsey$(s,t,n)$-graph for some $n$. There are finite number of Ramsey$(s,t)$-graphs for each $s$ and $t$~\cite{AS2008}, but finding all such graphs, or even determining the largest $n$ for which they exist, is a difficult problem. McKay's repository stores a large number of Ramsey$(s,t)$-graphs for different combinations of $s$ and $t$. 

McKay's repository stores all Ramsey$(3,4)$-graphs, all Ramsey$(3,5)$-graphs, all Ramsey$(3,6)$-graphs, all Ramsey$(3,7)$-graphs\footnote{Brinkmann et al.\ found $\num{1118436}$ graphs from this list \cite{BGP}}, all Ramsey$(3,8)$-graphs\footnote{Brinkmann and Goedgebeur found the full list in 2012 \cite{BGP}}, all maximal Ramsey$(3,9)$-graphs\footnote{The maximal Ramsey$(3,9)$-graph has 35 vertices and was found by Kalbfleisch in 1966~\cite{JK66}, but it took 47 years to prove its uniqueness \cite{JS}}, all Ramsey$(4,4)$-graphs, and all maximal Ramsey$(4,5)$-graphs.\ In 1995, McKay and Radziszowski proved that there are no Ramsey$(4,5)$-graphs with more than 24 vertices and found $\num{350904}$ of them with 24 vertices.\ The remainder were found in 2016 by McKay and Angeltveit. There are $\num{352366}$ altogether.

A \emph{Ramsey$(4,4;3)$-hypergraph} is a 3-uniform hypergraph (all hyperedges have size $3$) with this property: every set of $4$ vertices contains $1, 2$ or $3$ edges. The smallest order for which no such hypergraph exists is called the {\it hypergraph Ramsey number} R$(4,4;3)$. McKay computed all Ramsey$(4,4;3)$-hypergraphs up to $12$ vertices. 

\item{\bf Trees sorted by diameter} \\
McKay stored all possible trees up to $22$ vertices in a text file (though his program enables generation of all trees up to much higher order). The file contains all the trees of order $N$ and diameter $D$. There is one tree per line. The trees are given as an obvious list of edges, with vertices numbered from $0$. He also stored all homeomorphically irreducible trees up to $30$ vertices. They are also called series-reduced trees.\\

\item{\bf Digraphs}\\
{\tt digraph6} is a format used for storing directed graphs similar to the format used to store undirected graphs. McKay used it to store all the non-isomorphic tournaments\footnote{A {\em tournament} is a digraph obtained by assigning a direction for each edge in an undirected complete graph.} up to 10 vertices. A tournament of odd order $n$ is \emph{regular} if the out-degree of each vertex is $(n-1)/2$. A tournament of even order $n$ is {\em semi-regular} if the out-degree of each vertex is $n/2-1$ or $n/2$. McKay stores all regular and semi-regular tournaments of order up to 13.  

A regular tournament is \emph{doubly-regular} if each pair of vertices is jointly connected to exactly $(n-3)/4$ others~\cite{reidd}. The order of doubly-regular tournaments is $4n-1, n\in \mathbb{N}$. These tournaments are equivalent to skew Hadamard matrices~\cite{reidd}. McKay computed and stored these graphs up to 51 vertices, however, the list is incomplete for graphs of order $>27$. 

Let $G$ be a digraph and $v\in V(G)$. Let $N^+(v)=\{u\mid (u,v) \in E(G)\}$ and $N^-(v)=\{u\mid (v,u) \in E(G)\}$. A tournament is \emph{locally-transitive} if, for each vertex $v$, $N^+(v)$ and $N^-(v)$ are both transitive tournaments\footnote{A tournament is {\em transitive} if it is acyclic.}. McKay listed all the non-isomorphic locally-transitive tournaments up to $20$ vertices.  

McKay also computed acyclic digraphs up to 8 vertices. 
\end{itemize}
\begin{table}[htp]
\centering
\begin{tabular}{ |l|c|c| } 
\hline
Graph type & Order & Exhaustive list (Y/N)\\
\hline
\hline
Simple connected unlabelled& 1--11 & Y \\
\hline
Eulerian, Connected chordal  & 1--12, 1--13 & Y \\
\hline
Perfect&1--11&Y \\ \hline
Strongly regular  & refer~\cite{MK19,MS} & Y \\
\hline
Hypohamiltonian & 1--16& Y \\
\hline
Hypohamiltonian cubic& 1--26 & Y \\
\hline
Hypohamiltonian cubic, girth $\geq 5$ & 1--28 & Y \\
\hline
Hypohamiltonian cubic, girth $\geq 6$ & 1--30 & Y \\
\hline
Connected planar  & 1--11 & Y \\
\hline
Plane 5-regular simple connected & 1--36 & Y \\ \hline
Self-complementary & 1--17, 20& Y,N \\ \hline
Highly irregular  & 1--15& Y \\
\hline
Ramsey(3,4), Ramsey(3,5), Ramsey(3,6), & -- & Y \\
Ramsey(3,7), Ramsey(3,8), Ramsey(4,4)& &  \\ \hline
Maximal Ramsey(3,9), maximal Ramsey(4,5)& -- & Y \\
\hline
Tournament & 1--10 & Y \\ \hline
Regular tournament & 1--13 & Y \\ \hline
Semi-regular tournament & 1--13 & Y \\ \hline
Doubly-regular tournaments & 1--51 & N (The list is exhaustive \\ && with order up to 27) \\ \hline
Locally-transitive tournament&1--20 & Y \\ \hline
Acyclic digraph& 1--8 & Y \\ \hline
\end{tabular}
\caption{A summary of graphs from McKay's combinatorial data.}\label{mga}
\end{table}

All these graphs are presented in the form of static HTML pages and some of them can be downloaded as plain text files. McKay's repository is recognized for its sheer volume and variety. It expanded over the last three decades based on his research requirements. An outline of his repository is presented in Table~\ref{mga}. It is a very good resource of special classes of graphs. McKay's repository does not record graph parameters. Searching graphs of a special class (e.g., find a hypohamiltonian graph) is straightforward in McKay's combinatorial data, however it does not facilitate searching for graphs satisfying more than one property (e.g., find a graph which is hypohamiltonian and self-complementary). 
\subsubsection{Conder's combinatorial data (2002--present)}
Marston Conder maintains a collection of combinatorial group data, graphs and graph embeddings. In this section we will only discuss the graph-related data.  

Conder's collection of cubic graphs is arguably the largest collection of cubic graphs. Conder stores all cubic symmetric graphs up to isomorphism, on up to $\num{10000}$ vertices. One interesting thing to note here is all cubic symmetric graphs except the Petersen graph and the Coxeter graph\footnote{The \emph{Coxeter graph} is a 3-regular graph with 28 vertices and 42 edges~\cite{coex}.} have a Hamiltonian cycle. For each of these graphs, Conder computes their type (defined in~\cite{conda}), size of automorphism group, girth, diameter and bipartiteness. Conder also computed all symmetric graphs up to isomorphism, of order 2 to 30, with some information about their automorphism groups.
A regular graph is called \emph{semi-symmetric} if it is edge-transitive but not vertex-transitive (and then the automorphism group has two orbits on arcs). Conder computed all cubic semi-symmetric graphs up to isomorphism, on up to $\num{10000}$ vertices, listed by order, type (refer~\cite{conda}), girth, and diameter. 

Conder's contribution towards graph embeddings is also noteworthy. Conder gives a summary of the regular and chiral maps and maximum orders of group actions on compact Riemann surfaces of genus up to 301 and on compact non-orientable Klein surfaces of genus up to 302. For definitions and details on these data, refer to~\cite{MC2002}.  
    
Conder's repository focuses on graph embedding and symmetry data. Like McKay, Conder stored the graphs in static HTML pages. Conder's repository is a unique repository as it is one of the few repositories with a rich collection of geometric and combinatorial group-related graph data.

\subsubsection{Royle's combinatorial catalogue (2007--present)}
Gordon Royle's catalogue has a wide variety of parameters. It contains interesting data about other combinatorial objects like those used in finite geometry and design theory. In this survey we will focus on the graph-related data of the repository. We will briefly describe the scope of the repository and various graph data stored in the repository.

\begin{itemize}
\item {\bf Small graphs}\\
Royle's catalogue gives the total number of some specific types of graphs. They are listed in Table~\ref{t3}.
\begin{table}[htp]
\centering
\begin{tabular}{ |l|r| } 
\hline
\# of graphs of type & Up to order \\
\hline
\hline
Unlabelled simple& 16 \\
\hline
Unlabelled simple connected  & 16  \\
\hline
Multigraphs\tablefootnote{with number of edges varying from 1 to 14.} & 15\\
\hline
Connected bipartite\tablefootnote{This includes graphs with partitions $(1,14),(2,13),\cdots,(7,8)$. } & 15  \\
\hline
Trees & 20 \\
\hline
Class 1 and Class 2 \tablefootnote{As per Vizing's theorem.} & 9 \\
\hline
\end{tabular}
\caption{Royle's data on graph counting.}\label{t3}
\end{table}
  
Royle also computed some interesting parameters for the graphs.  This repository lists the chromatic numbers of all connected graphs on up to 11 vertices. A graph is said to be \emph{vertex-critical} if its chromatic number drops whenever a vertex is deleted. Royle lists all such graphs up to 11 vertices. The graphs are stored in gzip-compressed files containing those graphs in {\tt graph6} format. A graph is said to be \emph{edge-critical} if its chromatic number drops whenever an edge is deleted. Every edge-critical graph is necessarily vertex-critical. Royle lists all such graphs up to 12 vertices, he also  lists the number of all 4, 5, 6 and 7 edge-critical graphs up to 12 vertices, specifying the number of edges in each graph~\cite{GR2007}. Each  graph is stored in {\tt graph6} format.

\item{\bf Cubic graphs}\\
One feature of Royle's repository is that he not only stores the graphs but also computes (stores) some interesting parameters of the graphs. In the case of cubic graphs, Royle has listed all graphs with order up to 22 vertices with girth ranging from 2 to 14, he also listed all 3-connected cubic graphs (irrespective of girth) with up to 20 vertices. 

A \emph{snark} is defined to be a cyclically $4$-edge connected $3$-regular graph with chromatic index $4$ and girth $\ge 5$~\cite{snark}. The Petersen graph is the smallest snark, and Tutte conjectured that every snark has a Petersen graph minor~\cite{tu66}. This conjecture was proven in 2001 by Robertson, Sanders, Seymour, and Thomas, using an extension of the methods they used to reprove the four-colour theorem. Royle used Gunnar Brinkmann's cubic graph generation program to construct snarks of all orders up to 28, which are stored in the {\tt graph6} format.

We have discussed chromatic polynomials in Section~\ref{aag}. Royle used the tree representation to represent chromatic polynomials.
Given a family of polynomials, it is natural to explore their roots and to see whether these have any patterns. The complex roots of chromatic polynomials of graphs have interesting patterns, many of which are not explained yet. Read and Royle~\cite{RR91} investigated some of these patterns. They also calculated the chromatic roots of many small graphs. One of the longstanding questions was whether it was possible for the chromatic polynomial of a graph to have a root whose real part is negative~\cite{cp1,cp2}. Read and Royle discovered that certain cubic graphs of girth $\geq 5$ on $18$ vertices, girth $\geq 6$ on $20$ vertices, and girth $\geq 7$ on $26$ vertices have roots with quite significant negative real parts~\cite{RR91}. Royle's repository stores the chromatic polynomial of all connected cubic graphs up to 30 vertices with girth ranging from 3 to 8.

\item{\bf Transitive graphs}\\
Like McKay's repository, Royle's repository also contains transitive graphs. The data was prepared by McKay, Royle and  Alexander Hulpke. McKay generated transitive graphs up to 19 vertices~\cite{mtrans79}, subsequently Royle extended the data first up to 24 vertices~\cite{roylethesis} and then up to 26 vertices~\cite{mcroyle}. The current extension is done by Alexander Hulpke~\cite{ah30} who has constructed all the transitive groups of degree up to 30. Using these groups Royle performed a complete re-computation of the graphs. The re-computation confirmed that the original numbers computed by McKay and Royle were correct. 

Royle's repository stores all transitive graphs up to 31 vertices and he verified the correctness of the results up to 26 vertices. Royle also provided information on whether or not  the following properties hold for these transitive graphs: Cayley, Non-Cayley, Connected Transitive,	Connected Cayley, and Connected Non-Cayley hold on these graphs. 

\item{\bf Cayley graphs}\\
Let $G$ be a group, and let $S \subseteq G$ be a set of group elements such that the identity element $I$ is not in $S$. The \emph{Cayley graph} associated with $(G,S)$ is then defined as the directed graph having one vertex associated with each group element and directed edges $(g,h)$ whenever $gh^{-1}$ in $S$. The Cayley graph depends on the choice of the set $S$, and is connected if and only if $S$ generates $G$. Royle lists all the Cayley graphs on up to 31 vertices, but classified according to the group to which they belong.

The first group of each order is the cyclic group, then the remaining groups are ordered according to the lists in the book Group Tables by Thomas and Wood~\cite{shiva}. They give some descriptive names, which Royle also uses (with the exception of using $D(2n)$ for the dihedral group of order $2n$, rather than $D(n)$). All graphs are stored in the {\tt graph6} format.

\item {\bf Cubic transitive graphs}\\
The combinatorial catalogue of cubic vertex-transitive graphs (prepared by McKay and Royle) contains graphs of order up to 256 vertices (inclusive). The graphs are stored in {\tt graph6} format. Royle also specifies whether graphs are Cayley, non-Cayley and symmetric. However the list does not contain the information for some graphs with more than 96 vertices, if they are Cayley or non-Cayley.

\item {\bf Cubic cages and higher valency cages}\\
A $(k,g)${\it-cage} is a $k$-regular graph of girth $g$ with the fewest possible number of vertices.
Royle lists the currently known values for the sizes of a cubic cage. For certain small values of $g$ the cages themselves are all known, and Royle lists them explicitly in the repository. For larger values of $g$ Royle has given a range --- the lower value is either the trivial bound $n(3,g)$ or a bound by extensive computation. Royle lists some cubic $(3,g)$ graphs for $g \leq 32$. The cubic cages are stored in {\tt sparse6} format.

\item{\bf Cubic planar graphs}\\
Royle lists 3-connected cubic planar graphs with up to $20$ vertices. Royle constructed these by using Brinkmann and McKay's program {\tt plantri}~\cite{brmc}.

\item{\bf Cubic symmetric graphs (The Foster Census)}\\
A graph is called \emph{symmetric} if its automorphism group acts transitively on the set of arcs (directed edges) of the graph.

If the graph is cubic, then by Tutte's theorem~\cite{cututte,tuttesy} the automorphism group actually acts regularly on $s$-arcs for some value of $s$ between 1 and 5, and we say that a graph is {\em $s$-arc transitive} if the group acts regularly on $s$-arcs but not transitively on $(s+1)$-arcs. This repository lists the known cubic symmetric graphs with less than 1000 vertices. 

\item {\bf Strongly regular graphs}\\
Royle's repository stores strongly regular graphs up to 99 vertices and the graphs are stored in {\tt graph6} format. Royle also computed some useful parameters for each of these graphs. He computed the eigenvalues of each of these graphs with their multiplicity. 
\end{itemize}

\begin{table}[htp]
\centering
\begin{tabular}{ |l|c|c| } 
\hline
Graph type & Order & Exhaustive list (Y/N)\\ \hline
\hline

Vertex-critical& 1--11 & Y \\
\hline
Edge-critical& 1--12 & Y \\
\hline
Class-2 & 1--9 & Y \\ \hline
Cubic, diameter $\in [2,14]$  & 1--22 & N \\
\hline
3-connected cubic & 10--20 & Y \\ \hline
Snark & 1--28&Y \\ \hline
Transitive & 1--26 \tablefootnote{Royle has computed transitive graphs with up to 31 vertices, but correctness has been checked up to 26 vertices.}& Y \\ \hline
Cubic transitive &1--256& Y \\ \hline
Cayley & 1--31 & Y \\ \hline
Cubic cages ((3,3)-cage, (3,4)-cage,& -- & N \\ 
(3,5)-cage, $\cdots$ (3,22)-cage ) & & \\ \hline
3-connected cubic planar & 10--20 & Y \\ \hline
Cubic symmetric & 1--1000& N \tablefootnote{The list is complete up to order 768; for order in the range 770--798 it includes only Cayley graphs.} \\ \hline
Strongly regular graphs & 1--99 & N \\
 \hline
\end{tabular}
\caption{A summary of graphs from Royle's repository.}\label{rga}
\end{table}

Although Royle's repository contains fewer graphs than McKay's, it has  great variety and contains information about many graph parameters. A short summary of the graphs stored in Royle's repository is given in Table~\ref{rga}. There is some duplication in Royle's repository, and just like McKay's repository, the data is represented in static form.

\subsubsection{Other relevant repositories}
A graph $G$ is a {\em topological obstruction} for the torus if $G$ has minimum degree at least three, and $G$ does not embed on the torus but for all edges $e$ in $G$, the subgraph $G\setminus e$ embeds on the torus. A graph $G$ is a {\em minor-order obstruction} for the torus if $G$ is a topological obstruction for the torus and for all edges $e$ in $G$, the graph resulting from contracting $e$ embeds on the torus. Wendy Myrvold computed minor-order obstructions for the torus with up to 26 vertices \cite{ENS,ENM}. The graphs are stored as the upper triangular part of their adjacency matrix.

Edwin van Dam and Ted Spence worked on the classification of all regular graphs on at most 30 vertices that have four distinct eigenvalues~\cite{VD}. They classified the graphs in two parts: graphs for which all four eigenvalues are integral, and the case when there are just two integral eigenvalues. Spence also listed strongly regular graphs on at most 64 vertices. Spence contributed immensely to finding SRGs. Spence's findings with McKay are summarised in Section~\ref{bdm}. Along with Coolsaet and Degraer, Spence has computed the (45,12,3,3) strongly regular graphs~\cite{TED}. There are precisely 78 of these listed in the repository. 

Markus Meginger focuses on regular graphs. Meginger lists   simple connected $k$-regular graphs on $n$ vertices and girth at least $g$ with given parameters $n,k,g$ \cite {MM}. Meginger lists these graphs by using a computer program {\tt GENREG}. It not only computes the number of regular graphs for the chosen parameters but even constructs the desired graphs. The large cases with $k=3$ were solved by Gunnar Brinkmann, who implemented a very efficient algorithm for cubic graphs~\cite{BJ2013}. Meginger's collection of regular graphs are given in Table~\ref{mag}.
\begin{table}[htp]
\centering
\begin{tabular}{ |l|l|c| } 
\hline
Graph type & Order & Exhaustive list (Y/N)\\ \hline \hline
Connected 3-regular & 1--18 & Y \\ \hline
Connected 4-regular & 1--14 & Y \\ \hline
Connected 5-regular & 1--12 & Y \\ \hline
Connected 6-regular & 1--11 & Y \\ \hline
Connected 7-regular & 1--11  & Y \\ \hline
Connected 3-regular with girth $\geq 4$& 1--20 & Y \\ \hline
Connected 4-regular with girth $\geq 4$& 1--16 & Y \\ \hline
Connected 5-regular with girth $\geq 4$& 1--16 & Y \\ \hline
Connected 6-regular with girth $\geq 4$& 1--18 & Y \\ \hline
Connected 7-regular with girth $\geq 4$& 1--8  & Y \\ \hline
Connected 3-regular with girth $\geq 5$& 1--22 & Y \\ \hline
Connected 4-regular with girth $\geq 5$& 1--23 & Y \\ \hline
Connected 5-regular with girth $\geq 5$& 1--30 & Y \\ \hline
Connected 3-regular with girth $\geq 6$& 1--24 & Y \\ \hline
Connected 4-regular with girth $\geq 6$& 1--34 & Y \\ \hline
Connected 3-regular with girth $\geq 7$& 1--32 & Y \\ \hline
Connected 3-regular with girth $\geq 8$& 1--40 & Y \\ \hline
 \end{tabular}
 \caption{Meginger's repository.}\label{mag}
 \end{table}

Brinkmann's repository contains numbers of connected regular graphs with given number of vertices (up to 26) and degree (up to degree 7, for graphs of order 17 and of degree 3 for the rest). It also contains connected regular graphs with girth at least 4 for graphs of order 26, connected regular graphs with girth at least 5 for graphs of order 32, connected regular graphs with girth at least 6 for graphs of order 34, connected regular graphs with girth at least 7 for graphs of order 32, connected regular graphs with girth at least 8 for graphs of order 40 and connected bipartite regular graphs up to 32 vertices. Brinkmann's repository also contains connected planar regular graphs, along with connected planar regular graphs with girth at least 4 up to 26 vertices (degree 3), and connected planar regular graphs with girth at least 5 of order 20--28 with degree 3.

Primo\v z Poto\v cnik provides censuses of cubic and 4-regular graphs having different degrees of symmetry. Poto\v cnik's repository consists of the following graphs.
\begin{itemize}
\item {\bf Census of rotary maps}:  
For definitions of chiral, orientable maps and non-orientable maps, refer to~\cite{GT, bonningtonl}. Poto\v cnik's census stores all non-orientable maps up to 1500 edges and all orientable maps (both chiral and regular) up to 3000 edges. The previously known census of Conder contained the maps up to 1000 edges. 
\item {\bf Census of cubic vertex-transitive graphs}: A census (compiled by Pablo Spiga, Gabriel Verret and Primo\v z Poto\v cnik) of cubic vertex-transitive graphs on at most 1280 vertices is available at this link [\url{http://www.matapp.unimib.it/\~spiga/census.html}]. 
\item {\bf Census of 2-regular arc-transitive digraphs}: Spiga, Verret and Poto\v cnik have compiled a complete list of all connected arc-transitive digraphs on at most 1000 vertices \cite{PP}. As a byproduct, they have computed all connected 4-regular graphs with at most 1000 vertices that admit a half-arc-transitive action of a group of automorphisms. In particular, all 4-regular half-arc-transitive graphs are there.

\item{\bf Census of arc-transitive 4-regular graphs}: Spiga, Verret and Poto\v cnik, were able to construct a complete list of all 4-regular arc-transitive graphs on at most 640 vertices \cite{PPG,PSG}. The MAGMA code which generates the sequence of these graphs can be found at \url{https://www.fmf.uni-lj.si/\~potocnik/work_datoteke/Census4val-640.mgm}.

\item{\bf Census of 2-arc-transitive 4-regular graphs}: This repository contains a list of 2-arc-transitive 4-regular graphs on up to 2000 vertices (in MAGMA code). The list is complete on up to 727 vertices, but misses some 7-arc-transitive graphs that admit no $s$-arc-transitive group for $s$ less than 7 on more than 727 vertices, and also some 4-arc-transitive graphs that admit no $s$-arc-transitive group for $s$ less than 4 on more than 1157 vertices \cite{POT}.

Gary Haggard~\cite{ghag,ghag1,ghag2} computed and stored chromatic polynomials and the Tutte polynomials for some special classes of graphs. We summarise the contents of Haggard's repository in Table~\ref{hagre}.
\begin{table}[htp]
\centering
\begin{tabular}{ |l|c|c| } 
\hline
Graph type & Chromatic polynomial & Tutte polynomial\\ \hline \hline
9-cage-$k, k \in[1,18]$ & Y & N \\ \hline
Complete graphs of order $k$, $k \in [3,20]$ & Y & N \\ \hline
Complete graphs of order $k$, $k \in [6,25]$ & N & Y \\ \hline
\end{tabular}
\caption{A summary of Haggard's repository.}\label{hagre}
\end{table}

Suppose $G$ is a graph on $n$ vertices with diameter $d$. For any vertex $u$ and for any integer $i$ with $0\le i \le d$, let $G_i(u)$ denote the set of vertices at distance $i$ from $u$.
If $v\in G_i(u)$ and $w$ is a neighbour of $v$, then $w$ must be at distance $i-1$, $i$ or $i+1$ from $u$. Let $c_i, a_i  ~\hbox{and}~ b_i$ denote the number of vertices whose distances from $w$ are $i-1$, $i$ and $i+1$ respectively. $G$ is a {\it distance-regular graph} if and only if these parameters $c_i,a_i,b_i$ depend only on the distance $i$, and not on the choice of $u$ and $v$. R. Bailey, A. Jackson and C. Weir~\cite{bailey} developed a repository of distance-regular graphs. A summary of the repository is presented in Table~\ref{bolga}.\ Bailey et al.\ also maintains an index for all the named graphs present in the repository.
\begin{table}[htp]
\centering
\begin{tabular}{ |l|c| } 
\hline
Graph type & Exhaustive list (Y/N)\\ \hline \hline
Special families graphs & N \\ \hline
Unlabelled graph with up to order 1416 & N \\ \hline
 Distance-regular graphs with degree $k$, $k \in [3,13]$ & Y \\ \hline
Distance-regular graphs with diameter $k$, $k \in [2,10]$ & Y \\ \hline
Special named graphs & N \\ \hline

\end{tabular}
\caption{A summary of Bailey et al.'s repository.}\label{bolga}
\end{table}

All these repositories are rich in content and are useful resources in their own right, however they are not intended to be general graph repositories and do not cover many basic graph parameters. A summary of these repositories is given in Table~\ref{otga}. Moreover, these repositories are all static in nature, so they are not well suited to queries consisting of combinations of parameters.

\end{itemize}
\begin{table}[htp]
\centering
\begin{tabular}{ |l|l|c|c| } 
\hline
Repository name &Graph type & Order & Exh. list (Y/N)\\ \hline
Myrvold &Minor order & 1--26 & Y \\
&obstruction for torus&&\\
\hline
Spence &Strongly regular& refer~\cite{TED,MS} & N \\
\hline
Spence & Regular & 1--30 & N (refer to~\cite{VD}) \\ \hline 
Brinkman & Connected regular,  & 1--32 & N \\
 & girth $\ge k~ (k\in[4,7])$ & &  \\ \hline

Poto\v cnik & Rotary maps & -- \tablefootnote{up to 1500 edges for non-orientable maps and up to 3000 edges for orientable maps} & Y \\ \hline
Poto\v cnik et al.& 2-regular arc-transitive digraph& 1--1000 & Y \\ \hline
Poto\v cnik et al.& Arc-transitive 4-regular & 1--640 & Y \\ \hline
Poto\v cnik et al.& 2-arc-transitive 4-regular & 1--2000 & N \tablefootnote{Complete up to 727 vertices} \\ \hline
\end{tabular}
\caption{A summary of other useful repositories.}\label{otga}
\end{table}

\subsection{Interactive repositories}\label{irr}
The resources discussed above are static in nature. Although these resources are very useful they are not flexible, as the user has to manually refine the data if the requirements of the user is different to the data presented in the repository.

Data stored in these repositories are mostly in the form of adjacency matrix, adjacency list or McKay's {\tt graph6} and {\tt sparse6} format~\cite{g6}.\ To perform any tasks that involve these large data the user may need some built-in tools and computer(s). The advent of large scale computing clusters, GPUs, and high performance computers enables users to use large scale data in a more flexible and interactive way. In the following sections we list some relevant attempts at creating interactive graph repositories.

\subsubsection{Wolfram Alpha}
Wolfram Alpha (also styled WolframAlpha, and Wolfram|Alpha) is a computational knowledge engine or answer engine developed by Wolfram Alpha LLC, a subsidiary of Wolfram Research. It is an online service that answers factual queries directly by computing the answer from externally sourced ``curated data" \cite{AWA} rather than providing a list of documents or web pages that might contain the answer, as a search engine does.

Wolfram Alpha, which was released on May 18, 2009, is based on Wolfram's earlier flagship product Wolfram Mathematica, a computational platform or toolkit that encompasses computer algebra, symbolic and numerical computation, visualisation, and statistics capabilities. It contains definitions of more than 300 different types of undirected graphs, more than 20 different types of directed graphs and a  good collection of problems related to graphs. 

Although this is an informative repository, it is different from all other repositories we discuss. It has definitions and some theoretical information about graphs but it does not list or store graphs. It can be queried for gathering information about specific graphs, e.g., ``cliques" or ``Petersen graph". This repository does not filter graphs satisfying constraints like ``listing all graphs with max degree $k$". The number of graphs listed in Wolfram Alpha is smaller than the number of graphs listed in ``An Atlas of Graphs".

\subsubsection{Encyclopedia of Graphs (2012--present)}
Encyclopedia of Graphs (\url{http://atlas.gregas.eu/}) is an online encyclopedia of graph collections aiming to help researchers find and use data about various families of graphs. This repository  was created by a Slovenian company, Abelium d.o.o.. This repository allows graphs to be stored in any of these formats: {\tt graph6}, {\tt sparse6}, adjacency matrix, adjacency list and edge list. Encyclopedia of Graphs lists cubic symmetric graphs, edge-transitive graphs, vertex-transitive graphs and other classes of graphs, some of them acquired from other graph repositories listed in Section~\ref{err}. This repository contains the following graphs.
\begin{itemize}
\item {\bf Symmetric cubic graphs (The Foster Census)}: This repository stores the 796  cubic connected symmetric graphs with up to 2048 vertices.  
\item{\bf Edge-transitive 4-regular graphs}: The collection provides information about connected edge-transitive graphs of degree 4. The current edition has 793 graphs with up to 150 vertices, which turned out to be an incomplete list. The authors of this repository are aiming to expand the range up to 512 vertices.
\item{\bf Vertex-transitive graphs}: This page lists all the transitive graphs on up to 31 vertices. The data was prepared by  McKay and Royle~\cite{MK19,GR2007}, and Hulpke. The data in the current version are guaranteed to be correct  up to 26 vertices. The transitive groups on 24, 27, 28 and 30 vertices have not yet been checked. This list contains $\num{100661}$ graphs.

\item{\bf Hexagonal capping of symmetric cubic graphs} : The \emph {hexagonal capping} HC(G) of a graph $G$ has four vertices $\{u_0, v_0\},\{u_0, v_1\},\{u_1, v_0\},\{u_1, v_1\}$ for each edge $\{u,v\}$ of $G$, and each $\{u_i, v_j\}$ is joined to each $\{v_j, w_{1-i}\}$, where $u$ and $w$ are distinct neighbours of $v$ in $G$~\cite{hill}. This collection includes hexagonal cappings of cubic symmetric graphs up to 798 vertices. This list contains 284 graphs.

\item{\bf Line graphs of symmetric cubic graphs}: Given a graph $G$, its \emph{line graph $L(G)$} is a graph such that each vertex of $L(G)$ represents an edge of $G$; and two vertices of $L(G)$ are adjacent if and only if their corresponding edges share a common endpoint are incident in $G$. This repository includes line graphs of cubic symmetric graphs with up to 768 vertices and Cayley graphs in the range 770--998 vertices. 

\item{\bf Arc-transitive 4-regular graphs}: This collection contains a complete census of all connected arc-transitive 4-regular graphs of order at most 640. The census is a joint project by Poto{\v c}nik, Spiga, and Verret \cite{PSG}. They stored $\num{4820}$ graphs in this category.

\item{\bf Regular graphs}: This collection contains all connected regular graphs of girth 3 up to order 12, of girth 4 up to order 16, and of girth 5 up to order 23. The collection was prepared using the program {\tt geng} which is part of the {\tt nauty} software package of  McKay. This list contains $\num{140959}$ graphs.

\item{\bf Trees}: This collection contains  all $\num{522958}$ non-trivial trees with up to the order 19, as prepared by Royle \cite{GR2007}.

\item{\bf Vertex-transitive cubic graphs}: This collection contains a complete census of all $\num{111360}$ connected vertex-transitive cubic graphs of order at most $\num{1280}$. The census is a joint project by Poto{\v c}nik, Spiga, and Verret \cite{PSG}.\ The properties of graphs were collected from the DiscreteZOO library \cite{DZ}. 

\item{\bf Highly irregular graphs}: A connected graph is \emph{highly irregular} if the neighbours of each vertex have distinct degrees. Such graphs exist on all orders except 3, 5 and 7. This collection lists all $\num{21869}$ highly irregular graphs with up to the order 15 and was provided by  McKay \cite{MK19}. 
 
\item{\bf Snarks}: This lists all $\num{153863}$ snarks with up to the order 30 and was provided by The House of Graphs \cite{BJ2013}.
 
\item{\bf Cubic graphs}: This collection contains a complete census of the $\num{556471}$ connected cubic graphs of order at most 20. The data was prepared by Royle \cite{GR2007} and later extended at the House of Graphs \cite{BJ2013}. 

\item{\bf Strongly regular graphs}: We defined strongly regular graphs in the Section~\ref{bdm}.  A strongly regular graph is called {\it primitive} if both the graph and its complement are connected. The census lists $\num{43679}$ primitive strongly regular graphs with order up to 40; however, the complete classification is still open for SRG(37,18,8,9). The census was provided by  McKay and Spence~\cite{MS}.

\item{\bf Networks}: Network theory is the study of graphs as  representations of either symmetric or asymmetric relations between discrete real-world objects. It studies technological networks (the internet, power grids, telephone networks, transportation networks), social networks (social graphs, affiliation networks), information networks (World Wide Web, citation graphs, patent networks), biological networks (biochemical networks, neural networks, food webs), and many more. Graphs provide a structural model that makes it possible to analyse and understand the ways in which many separate entities act together. This (expanding) collection houses some of the ``interesting" networks that are used in research. This category contains 51 graphs of order ranging from $3$ to $23219$. Most of these graphs are of order ranging from $30$ to $40$.

\item{\bf Edge-critical graphs}: The collection was calculated by Royle~\cite{GR2007} and lists the $\num{185844}$ edge-critical graphs with from 4 to 12 vertices.
\item{\bf Fullerenes}: A \emph{fullerene} is a cubic planar graph having all faces 5- or 6-cycles~\cite{grun}. Fullerenes are planar and hence polyhedral, and every fullerene has exactly twelve 5-cycles. They acquired this data from the House of Graphs~\cite{BJ2013} and listed all $\num{467927}$ fullerenes on up to 90 vertices. 

\item{\bf Maximal triangle-free graphs}: In order to determine properties of all triangle-free graphs, it often suffices to investigate maximal triangle-free graphs. These are triangle-free graphs for which the insertion of any further edge would create a triangle. They acquired this data from the House of Graphs \cite{BJ2013} and listed all the $\num{197396}$ maximal triangle-free graphs on up to 17 vertices. 

\item{\bf Planar graphs}: This collection was calculated by  McKay~\cite{MK19} and he lists all $\num{78633}$ non-isomorphic connected planar graphs on up to 9 vertices. 

\item{\bf Vertex-critical graphs}: The collection was calculated by Royle~\cite{GR2007} and lists the $\num{359787}$ vertex-critical graphs with from 4 to 11 vertices.
\end{itemize}

One can search the site by using graph names, Universal Graph Identifier (UGI), defined by the authors of the repository, or by collection names. UGI is a string that uniquely identifies a graph and can be used to directly access its properties page. One can use filters to obtain specific graphs of interest, e.g.\ ``bipartite = true, minimum degree $> 3$". After the relevant graphs and their properties have been selected, they can be exported. Alternatively, one can download the data of a specific graph. Its main focus is to list families of graphs. A summary of the graphs stored in this repository is listed in Table~\ref{ega}. Since most of its data has been collected from other sources, the correctness of the data is subject to the correctness of the various repositories used as sources.

\begin{table}[htp]
\centering
\begin{tabular}{ |l|l|l| } 
\hline
Graph type & Order & Exhaustive list (Y/N)\\ \hline
\hline
Connected Cubic & 1--20 & $\mbox{Y}^{*,\#}$ \\ \hline
Symmetric cubic & 1--2048 & $\mbox{Y}^{*, \%}$ \\ \hline
Transitive & 1--31 & $\mbox{Y}^{*,+}$ \\ \hline
Vertex-transitive cubic & 1--1280 & $\mbox{Y}^{\&}$ \\ \hline
Edge-transitive 4-regular & 1--150 & N \\ \hline
Arc-transitive 4-regular & 1--640 &  $\mbox{Y}^{\&}$ \\ \hline
Line graphs of symmetric cubic graphs & 1--768 & Y \\ \hline

Connected regular, girth 3 & 1--12 & $\mbox{Y}^{+}$ \\ \hline
Connected regular, girth 4 & 1--16 & $\mbox{Y}^{+}$ \\ \hline
Connected regular, girth 5& 1--23 & $\mbox{Y}^{+}$ \\ \hline
Trees & 1--19 & $\mbox{Y}^{*}$ \\ \hline
Maximal triangle-free graphs & 1--17 & $\mbox{Y}^{\#}$ \\ \hline
Planar  & 1--9 & $\mbox{Y}^{+}$\\
\hline
Vertex-critical & 4--11 & $\mbox{Y}^{*}$ \\ \hline
Edge-critical & 4--12 & $\mbox{Y}^{*}$ \\
\hline
Strongly regular & 1--40 & $\mbox{N}^{+}$ \cite{MS} \\ \hline
Highly irregular graph & 1--15 & $\mbox{Y}^{+}$ \\ \hline 
Snark & 1--30 & $\mbox{Y}^{\#}$ \\
\hline
Fullerenes & 1--90 & $\mbox{Y}^{\#}$ \\ \hline
Network & -- & N \\ \hline

\end{tabular}
\caption{A summary of graphs from Encyclopedia of Graphs.\ +: data from McKay's repository; *: data from Royle's repository; \%: data from Conder's repository; \#: data from house of graphs; \&: data from Poto\v cnik's repository.} \label{ega}
\end{table}

\subsubsection{House of graphs (2013--present)}
A recent attempt to build an interactive repository ``A House of Graphs" (\url{https://hog.grinvin.org/}) is by Brinkmann, Coolsaet, Goedgebeur and M\'elot in 2013~\cite{BJ2013}. A House of Graphs provides a searchable and downloadable graph database. Another functionality of the House of Graphs is a list of graphs which have been used as counterexamples. The authors of the database call these graphs ``interesting/relevant" and they suggest that, for a new theorem or conjecture, the chance of finding a counterexample is higher among the ``interesting" graphs. The authors acknowledge the difficulty in creating a database that contains all graphs up to a specific order, therefore they chose to create a database based on a paraphrase of Orwell's famous words: {\em all graphs are interesting, but some graphs are more interesting than others}~\cite{BJ2013}.

Some graphs (e.g., the Petersen graph or the Heawood $(3,6)$-cage on 14 vertices) or graph classes (e.g., snarks) appear repeatedly in the literature as counterexamples. 
In order to construct a rich source of possible counterexamples, they inserted 1570 graphs into the database at its inception. These graphs are either counterexamples to known conjectures or extremal graphs. The authors explicitly consider extremal graphs ``interesting".  A large proportion of the 1570 graphs are extremal graphs found by GraPHedron \cite{HM}. 

The idea behind creating this repository is that if one wants to test a conjecture on a list of graphs, the ideal case would be if one could restrict the tests to graphs that are ``interesting" or ``relevant" for this conjecture~\cite{BJ2013}. Here, the meaning of ``interesting" and ``relevant" is vague, but this already shows how much the question of whether a graph is interesting or not depends on the question that one wants to study. As per the authors, if a specific graph has sufficient properties to distinguish it in some way from the huge mass of other graphs, then the graph is interesting in some respect \textemdash~ and of course the database allows one to add that graph and also offers the possibility of saying for which invariants the graph is especially interesting. The House of Graphs offers (among others) the following lists:
\begin{itemize}
\item All graphs registered as interesting in the database. These interesting graphs plays a special role in this database. For graphs in this list, a lot of invariants (like the chromatic number, chromatic index, the clique number, the diameter, independence number, average degree, smallest eigenvalue, second largest eigenvalue, genus, etc.) and also embeddings (drawings) are precomputed and stored.
\item All snarks with girth at least 4 up to 34 vertices and with girth at least 5 up to 36 vertices.
\item All IPR-fullerenes\footnote{The {\it face-distance} between two pentagons is the distance between the corresponding vertices of degree 5 in the dual graph. We refer to the least face-distance between pentagons of a fullerene as the {\it pentagon separation} of the fullerene, denoted by $d$. Note that $d=1$ gives the set of all fullerenes and $d=2$ gives the set of all {\it IPR fullerenes}.}~\cite{ful} up to 400 vertices.
\item Complete lists of regular graphs for various combinations of degree, vertex number and girth.
\item Vertex-transitive graphs.
\item Some classes of planar graphs.
\end{itemize}

Some of these lists are physically situated on the same server as the website itself, but others are just links to other people's websites, like those of McKay~\cite{MK19}, Royle~\cite{GR2007}, Spence~\cite{TEDS} and Meringer~\cite{MM}.

Although this is not a static repository of graphs and it provides the option of querying the database with a combination of multiple parameters, it does not list all graphs up to some order. Graphs like Petersen's graph are obvious inclusions for this repository  but among special classes of graphs like trees, etc., very few will satisfy the ``interesting" criterion. As a result one cannot get a comprehensive list of graphs up to some order as an outcome of queries. We list the graphs of this repository in Table~\ref{hoga} and the parameters in Table~\ref{poga}. The repository allows users to add graphs to the database that they themselves consider ``interesting". This itself is an interesting way of building a searchable graph database.

\begin{table}[htp]
\centering
\begin{tabular}{ |l|l|l| } 
\hline 
Graph type & Order & Exhaustive\\
& & list (Y/N) \\ \hline
\hline
Snarks with girth $\geq 4$ & 1--34 & Y \\
\hline
Snarks with girth $\geq 5$ & 1--36 & Y \\
\hline
IPR-fullerenes & 1--160 & Y\\ \hline
Uniquely Hamilton with girth $\geq k$ ($k \in [3,5]$)& 1--12 &Y \\ \hline
Planar uniquely Hamilton with girth $\geq k$ ($k \in [3,5]$)& 1--12 &Y \\ \hline
Triangle-free $k$-chromatic   & 1--15, & Y \\ 
($k \in [4,5]$)& 1--23 & \\ \hline
Cubic  & 1--22 & $\mbox{N}^*$ \tablefootnote{All cubic graphs with diameter ranging from 2 to 14.} \\ \hline
Transitive & 1--26& $\mbox{Y}^*$ \\ \hline
Cubic transitive & 1--256 & $\mbox{Y}^*$ \\ \hline
Connected Cubic & 1--24 & $\mbox{Y}^{*,\$}$ \\ \hline
Connected regular (same as Table~\ref{mag}) & -- & $\mbox{Y}^{\$,+}$ \\ \hline
Strongly regular & -- & $\mbox{N}^{+,\#}$ \\ \hline
\end{tabular}
\caption{A summary of graphs from House of Graphs.\ +: data from McKay's repository; *: data from Royle's repository; \#: data from Spence's repository; \$: data from Meginger's repository.} \label{hoga}
\end{table}

\begin{table}[htp]
\centering
\begin{tabular}{ |l| } 
\hline
List of parameters\\ \hline \hline
Algebraic connectivity, average degree, circumference \\ \hline
Chromatic number, chromatic index, clique number \\ \hline
Density, diameter, edge connectivity \\ \hline
Genus, girth, longest induced path \\ \hline
Longest induced cycle, matching number, minimum independent set \\ \hline
Number of components, number of triangles, radius \\ \hline
Second largest eigenvalue, smallest eigenvalue, vertex connectivity\\ \hline
\end{tabular}
\caption{Parameters in House of Graphs.}\label{poga}
\end{table}

\subsubsection{Hoppe and Petrone's collection (2014)}
In 2014, Hoppe and Petrone~\cite{HOP} exhaustively enumerated all simple, connected graphs of order $\leq 10$ using {\tt nauty}~\cite{MK1984} and have computed the independence number, automorphism group size, chromatic number, girth, diameter and various properties like Hamiltonicity and Eulerianness over this set. Integer sequences were constructed from these invariants and checked against the Online Encyclopedia of Integer Sequences (OEIS). However they used brute force methods using Networkx~\cite{HS2008}, graph-tools~\cite{TP2014} and PuLP~\cite{MI2011} to compute the parameters. They presented the data in static form and stored it in a database. However this system is not interactive.

\subsubsection{Discrete ZOO}
Another notable repository is the Discrete ZOO  (\url{https://discretezoo.xyz/}). This repository hosts $\num{212269}$ graphs~\cite{DZ}. This repository mainly contains vertex transitive graphs (up to 31 vertices), cubic vertex transitive graphs (up to 1280 vertices) and cubic arc transitive graphs (up to 2048 vertices). One can filter its queries based on the graph parameters listed in Table~\ref{filter}. 
\begin{table}[htp]
\centering
\begin{tabular}{ |l|c| } 
\hline
{\bf Parameter} & {\bf Type}\\ \hline
\hline
Bipartite & Boolean \\
\hline
Cayley & Boolean \\
\hline
Clique number & Numeric \\
\hline
Degree & Numeric \\
\hline
Diameter & Numeric \\
\hline
Distance regular & Boolean\\
\hline
Distance transitive & Boolean \\
\hline
Edge transitive & Boolean \\
\hline
Girth & Numeric \\
\hline
Moebius ladder & Boolean \\
\hline
Strongly regular & Boolean \\
\hline
Triangles count & Numeric \\
\hline
\end{tabular}
\caption{Parameter filter used in Discrete ZOO.}\label{filter}
\end{table}
\subsubsection{An online atlas of graphs (2010--present)}\label{NBS}
The repositories by McKay~\cite{MK19}, Royle~\cite{GR2007}, Conder~\cite{MC2000} are rich in content and are considered very valuable resources for research. Any query on these data requires downloading the data first followed by manual compilation. The repositories are not designed for running queries consisting of multiple parameters on these data. The House of Graphs on the other hand is extremely efficient for handling complex queries but the data stored in this repository is sparse and incomplete. Although this is a good repository for some conjecture verification work, other problems may require more comprehensive data. 

Another approach is to give comprehensive information on queries (conjectures), i.e., {\it if} the system can answer precisely whether the conjecture holds for all graphs of order $\le k$, then this will at least become a lower bound on the order of graphs for which the conjecture holds, {\it else} the system will provide a counterexample for the conjecture.
This prompted Paul Bonnington and Graham Farr to propose a repository of graphs which has the efficiency of handling complex queries (like House of Graphs) and completeness of data (like the repositories of McKay, Royle etc.).

In 2009 Nick Barnes~\cite{NICK2010} built a prototype of an Online Graph Atlas (OLGA) with the following invariants: degree sequences, connected components, girth, radius and diameter, independence number, clique number, vertex cover number, domination number, circumference, length of the longest path, size of the maximum matching. Barnes used the Monash Grid cluster on a quad-core 2.5Ghz Intel Xeon L5420 processor,
with 16GB of RAM, and a MySQL server hosted by Monash University Information Technology Services to develop the prototype. All of these invariants that Barnes used are easily computable or they follow a recursive rule, but there are some parameters, like the genus of a graph, which do not follow a recursive structure.
C. Paul Bonnington and Graham Farr developed an approach to computing parameters using recursive lower and upper bounds. Man Son Sio (a BSE Honours student from Clayton School of IT, Monash University) extended Nick Barnes's system to include two new parameters: genus, as an example of a parameter that is NP-hard and has no known simple recursive rule; and the chromatic polynomial, as an example of a parameter that is a polynomial rather than a
single number. Sio also improved the database query times, which were a problem in the first prototype, though there is scope of further improvement there~\cite{SIO2010}.

The initial versions of OLGA used McKay's {\tt nauty}~\cite{MK19} as the backbone to generate the graphs. The advantage of using {\tt nauty} is the ease of use but it makes the system dependent on {\tt nauty}. In 2015, the author (under the supervision of Farr, Bonnington and Morgan) changed the basic architecture of OLGA. This work uses the Schreier-Sims algorithm~\cite{SC1971} for isomorphism checking and can generate all graphs up to 12 vertices. We took advantage of parallel computing and cloud computing. The current version of OLGA contains  degree sequence, connected components, girth, radius, and diameter, independence number, clique number, vertex cover number, domination number, circumference, length of the longest path, size of the maximum matching, vertex connectivity, edge connectivity, chromatic number, chromatic index, treewidth, Tutte polynomial (with order up to 7 vertices), automorphism group, genus, and eigenvalues of all graphs with order up to $10$ vertices. We also introduced a new graph parameter, the most frequent connected induced subgraph (MFCIS). The design of OLGA is scalable, however storage is a bottleneck for OLGA.

%
%
%
\section {In a nutshell}\label{nsh}
\begin{figure}
\center
 \includegraphics[width=\linewidth]{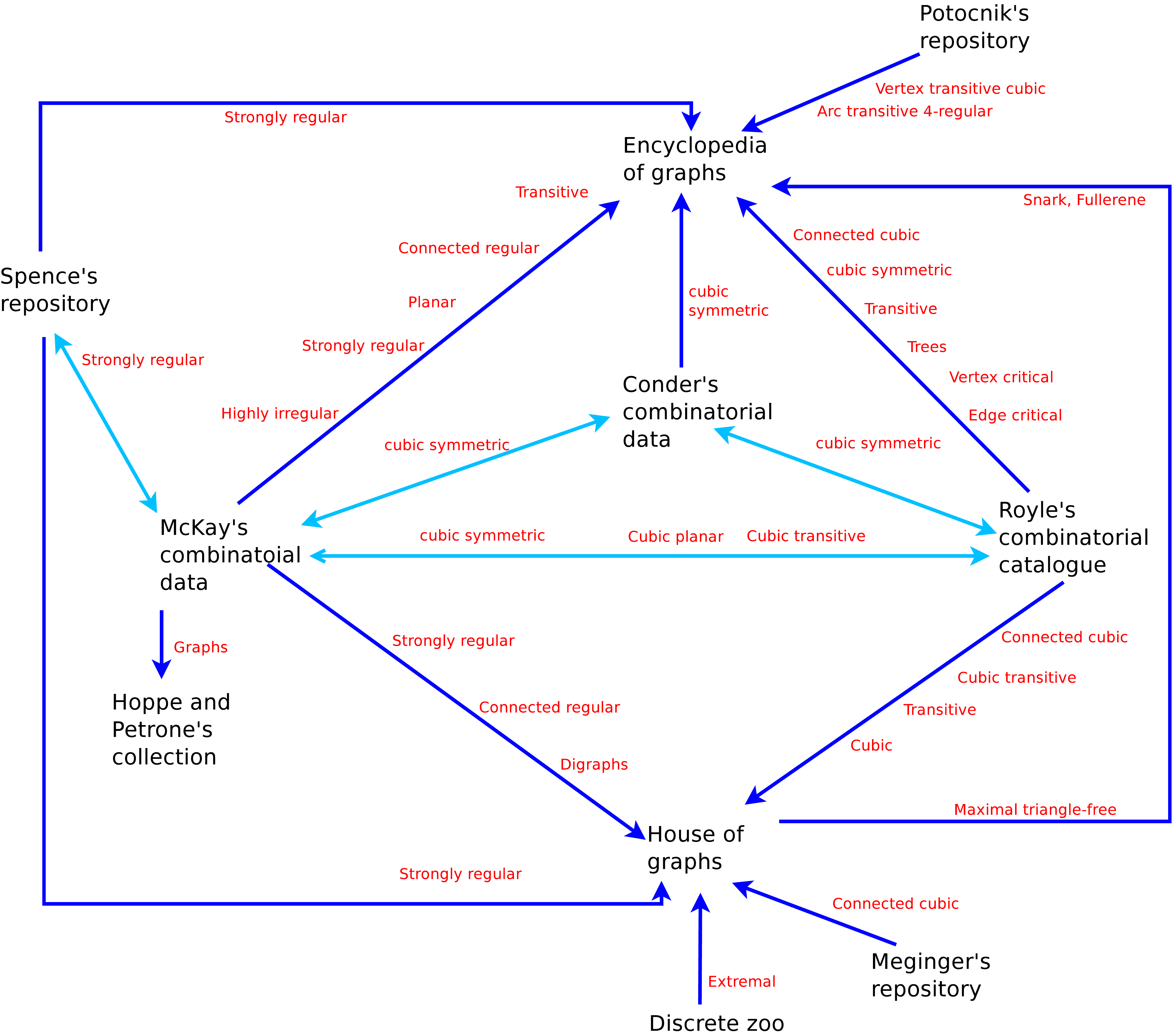}
\caption{Data flow between graph repositories. }
 \label{intd}
\end{figure}

In this section we summarise the contents of each repository in Table~\ref{sum}. In Figure~\ref{intd}, we illustrate the high level dependencies between the repositories using a simple directed labelled graph, with repositories as vertices and dependencies between repositories as edges. If repository ``B" depends on repository ``A" for some specific data ``info", we depict this by a directed edge from ``A" to ``B" with the label ``info" on it. If data is generated in collaboration between two repositories, we use a bidirectional edge between them.   

\begin{table}[htp]
\centering
\begin{tabular}{ |l|c| } 
\hline
{\bf Graphs and parameters} & {\bf Repositories}\\ \hline
Unlabelled graphs with up to order 10 & McKay, Royle \\ \hline
Strongly regular graphs & McKay, Spence, Royle, EoG \\ \hline
Ramsey graphs & McKay\\ \hline
Trees & McKay, Royle \\ \hline
Digraphs & McKay, HoG\\ \hline
Cubic graphs & Conder, Royle \\ \hline
Cubic symmetric graphs & Conder, Royle \\ \hline
Symmetric graphs & Conder \\ \hline
Snarks & Royle, HoG, EoG\\ \hline
Transitive graphs & Royle, Poto\v cnik, EoG \\ \hline
Cubic planar graphs & Royle, McKay \\ \hline
Fullerenes & House of Graphs, EoG\\ \hline
Extremal graphs & HoG, Discrete zoo \\ \hline 
Edge-critical graphs & Royle, EoG\\ \hline
Vertex-critical graphs & Royle, EoG\\ \hline

\end{tabular}
\caption{Summary of graph repositories. EoG = Encyclopedia of Graphs, HoG = House of Graphs.}\label{sum}
\end{table}
\section{Conclusion}\label{conc}
All these graph repositories are created to help researchers to get readily available graph data. If the need is to find a specific family of graphs like regular graphs, SRGs, Ramsey graphs etc., one can use the static repositories of McKay~\cite{MK19}, Royle~\cite{GR2007} and Conder~\cite{MC2000}. To get an overall picture and comprehensive information about a graph-related problem on graphs with up to 11 vertices, one can use OLGA. To verify theorems and conjectures on interesting graphs, A House of Graphs~\cite{BJ2013} can be used. None of these repositories guarantee a complete set of information about graphs or graph parameters, but they still offer useful information to understand many problems better and save a lot of time for researchers.  

\section*{Acknowledgement} The authors would like to thank Charles Colbourn, Marston Conder, and Brendan McKay for their comments on earlier versions of this paper.
\bibliographystyle{plain} 
\bibliography{myBibliography}

\begin{thebibliography}{10}

\bibitem{AS2008}
N.~Alon and J.~H. Spencer.
\newblock {\em The Probabilistic Method, 3rd edition}.
\newblock Wiley, 2008.

\bibitem{AR2009}
S.~Arora and B.~Barak.
\newblock {\em Computational Complexity: A Modern Approach}.
\newblock Cambridge University Press, 2009.

\bibitem{bailey}
R.~Bailey, A.~Jackson, and C.~Weir.
\newblock Distanceregular.org.
\newblock \\\url{http://www.distanceregular.org/index.html}.
\newblock Accessed: 2019-02-03.

\bibitem{bak1967}
G.~A. Baker, H.~E. Gilbert, J.~Eve, and G.~S. Rushbrooke.
\newblock {\em A Data Compendium of Linear Graphs with Application to the
  Heisenberg Model}.
\newblock BNL (Series). Brookhaven National Laboratory, 1967.

\bibitem{dewdney}
H.~H. Baker, A.~K. Dewdney, and A.~L. Szilard.
\newblock Generating the nine-point graphs.
\newblock {\em Mathematics of Computation}, 28(127):833--838, 1974.

\bibitem{NICK2010}
N.~Barnes.
\newblock {\em Towards an Online Graph Atlas}.
\newblock BCompSc Hons dissertation. Clayton School of Information Technology,
  Monash University, 2009.

\bibitem{DZ}
K.~Ber{\v{c}}i{\v{c}} and J.~Vidali.
\newblock Discrete {ZOO}: Towards a fingerprint database of discrete objects.
\newblock In {\em 6th International Congress on Mathematical Software}, volume
  10931 of {\em Lecture Notes in Computer Science}, pages 36--44, Editors: J.
  H. Davenport and M. Kauers and G. Labahn and J. Urban. South Bend, IN, USA,
  2018.

\bibitem{bonningtonl}
C.~P. Bonnington and C.~H.~C. Little.
\newblock {\em The Foundations of Topological Graph Theory}.
\newblock Springer, New York, 2012.

\bibitem{RAJ63}
R.~Bose.
\newblock Strongly regular graphs, partial geometries and partially balanced
  designs.
\newblock {\em Pacific Journal of Mathematics}, 13(2):389--419, 1963.

\bibitem{cage}
G.~Brinkmann, O.~D. Friedrichs, S.~Lisken, A.~Peeters, and N.~V. Cleemput.
\newblock {C}a{G}e - a virtual environment for studying some special classes of
  plane graphs - an update.
\newblock {\em MATCH Communications in Mathematical and in Computer Chemistry},
  63(3):533--552, 2010.

\bibitem{BJ2013}
G.~Brinkmann, J.~Goedgebeur, H.~M\'elot, and K.~Coolsaet.
\newblock House of {G}raphs: a database of interesting graphs.
\newblock {\em Discrete Applied Mathematics}, 161:311--314, 2013.

\bibitem{BGP}
G.~Brinkmann, J.~Goedgebeur, and J.~Schlage-Puchta.
\newblock Ramsey numbers ${R}({K}_3,{G})$ for graphs of order 10. arxiv
  preprint arxiv:1208.0501.
\newblock 2012.

\bibitem{brmc}
G.~Brinkmann and B.~D. McKay.
\newblock The program plantri.
\newblock \url{https://users.cecs.anu.edu.au/\~bdm/plantri/}, 2007.
\newblock Accessed: 2019-01-19.

\bibitem{TEDS}
A.~E. Brouwer and E.~Spence.
\newblock Cospectral graphs on 12 vertices.
\newblock {\em The Electronic Journal of Combinatorics}, 16(1):N20, 2009.

\bibitem{burr}
S.~A. Burr.
\newblock Diagonal {R}amsey numbers for small graphs.
\newblock {\em Journal of Graph Theory}, 7(1):57--69, 1983.

\bibitem{busmaker}
F.~C. Bussemaker, S.~Cobeljic, D.~M. Cvetkovic, and J.~J. Seidel.
\newblock Computer investigation of cubic graphs.
\newblock {\em EUT report. WSK, Dept. of Mathematics and Computing Science},
  76-WSK-01, 1976.

\bibitem{JA}
J.~D. Carpinelli and A.~Y. Oruc.
\newblock Applications of {M}atching and {E}dge-coloring {A}lgorithms to
  {R}outing in {C}los {N}etworks.
\newblock {\em Networks}, 24(6):319--326, 1994.

\bibitem{cp17}
N.~Cohen and D.~V. Pasechnik.
\newblock Implementing {B}rouwer's database of strongly regular graphs.
\newblock {\em Designs, Codes and Cryptography}, 84(1-2):223--235, 2017.

\bibitem{MC2002}
M.~Conder.
\newblock Combinatorial data.
\newblock \url{https://www.math.auckland.ac.nz/\~conder/}, 2002.
\newblock Accessed: 2016-01-03.

\bibitem{conda}
M.~Conder.
\newblock Trivalent (cubic) symmetric graphs on up to 10000 vertices.
\newblock
  \url{https://www.math.auckland.ac.nz/\~conder/symmcubic10000list.txt}, 2011.
\newblock Accessed: 2018-10-23.

\bibitem{MC2000}
M.~Conder and P.~Dobcs\'anyi.
\newblock Trivalent symmetric graphs on up to 768 vertices.
\newblock {\em Journal of Combinatorial Mathematics and Combinatorial
  Computing}, 40:41--63, 2002.

\bibitem{TED}
K.~Coolsaet, J.~Degraer, and E.~Spence.
\newblock The strongly regular $(45, 12, 3, 3) $ graphs.
\newblock {\em The Electronic Journal of Combinatorics}, 13(1):32, 2006.

\bibitem{VD}
E.~R.~V. Darn and E.~Spence.
\newblock Small regular graphs with four eigenvalues.
\newblock {\em Discrete Mathematics}, 189:233--257, 1998.

\bibitem{RD2010}
R.~Diestel.
\newblock {\em Graph Theory, 3rd edition}.
\newblock Springer, New York, 2010.

\bibitem{CD1960}
C.~Domb.
\newblock On the theory of cooperative phenomena in crystals.
\newblock {\em Advances in Physics}, 9(34--35), 1960.

\bibitem{BD}
R.~G. Downey and M.~R. Fellows.
\newblock {\em Parameterized Complexity}.
\newblock Springer-Verlag, 1999.

\bibitem{FB1988}
R.~M. Foster and I.~Z. Bouwer.
\newblock {\em The Foster Census : R.M. Foster's census of connected symmetric
  trivalent graphs}.
\newblock Winnipeg, Canada : Charles Babbage Research Centre, 1988.

\bibitem{frazer}
R.~J. Frazer.
\newblock {\em Graduate course project}.
\newblock unpublished, Department of Combinatorics and Optimization, University
  of Waterloo, 1973.

\bibitem{GJ1979}
M.~R. Garey and D.~S. Johnson.
\newblock {\em Computers and Intractability: A Guide to the Theory of
  NP-Completeness}.
\newblock W. H. Freeman \& Co., San Francisco, 1979.

\bibitem{ful}
J.~Goedgebeur and B.~D. McKay.
\newblock Fullerenes with distant pentagons.
\newblock {\em arXiv preprint arXiv:1508.02878}, cs.DS/1508.02878, 2015.
\newblock \url{https://arxiv.org/abs/cs.DS/1508.02878}.

\bibitem{JS}
J.~Goedgebeur and S.~P. Radziszowski.
\newblock New computational upper bounds for {R}amsey numbers ${R}(3,k)$.
\newblock {\em The Electronic Journal of Combinatorics}, 20(1):P30, 2013.

\bibitem{DW}
G.~Grimmett and C.~McDiarmid, editors.
\newblock {\em Combinatorics, {C}omplexity, and {C}hance: {A} {T}ribute to
  {D}ominic {W}elsh}.
\newblock Oxford Lecture Series in Mathematics and Its Applications 34, Oxford
  University Press, 2007.

\bibitem{GT}
J.~L. Gross and T.~W. Tucker.
\newblock {\em Topological Graph Theory}.
\newblock Dover Publications, Inc., New York, 1987.

\bibitem{grun}
B.~Gr{\"u}nbaum and T.~S. Motzkin.
\newblock The number of hexagons and the simplicity of geodesics on certain
  polyhedra.
\newblock {\em Canadian Journal of Mathematics}, 15:744--751, 1963.

\bibitem{HS2008}
A.~A. Hagberg, D.~A. Schult, and P.~J. Swart.
\newblock In {\em Exploring network structure, dynamics, and function using
  networkx}, pages 11--15, Editors: G. Varoquaux, J. Millman, T. Vaught.
  Pasadena, CA, 2008.

\bibitem{ghag1}
G.~Haggard.
\newblock Chromatic polynomials of 9 cages.
\newblock \\\url{http://www.eg.bucknell.edu/\~graphs/9cages.htm}.
\newblock Accessed: 2019-02-03.

\bibitem{ghag2}
G.~Haggard.
\newblock Chromatic polynomials of complete graphs.
\newblock \\\url{http://www.eg.bucknell.edu/\~graphs/complete.htm}.
\newblock Accessed: 2019-02-03.

\bibitem{ghag}
G.~Haggard.
\newblock Tutte polynomials of complete graphs.
\newblock \\\url{http://www.eg.bucknell.edu/\~graphs/tutte.htm}.
\newblock Accessed: 2019-02-03.

\bibitem{HP1973}
F.~Harary and E.~M. Palmer.
\newblock {\em Graphical Enumeration}.
\newblock Academic Press, 1973.

\bibitem{heap}
B.~R. Heap.
\newblock The production of graphs by computer.
\newblock In {\em Graph Theory and Computing}, pages 47--62. Editors: R. C.
  Read. Academic Press, Cambridge, Massachusetts, 1972.

\bibitem{hill}
A.~Hill and S.~Wilson.
\newblock Four constructions of highly symmetric tetravalent graphs.
\newblock {\em Journal of Graph Theory}, 71(3):229--244, 2012.

\bibitem{HOP}
T.~Hoppe and A.~Petrone.
\newblock Integer sequence discovery from small graphs.
\newblock {\em Discrete Applied Mathematics}, 201(C):172--181, 2016.

\bibitem{ah30}
A.~Hulpke.
\newblock Constructing transitive permutation groups.
\newblock {\em Journal of Symbolic Computation}, 39(1):1--30, 2005.

\bibitem{PH1964}
P.~Hutchinson.
\newblock Diagram expressions useful in the theory of fluids.
\newblock Technical report, United Kingdom Atomic Energy Authority (Research
  Group), 1964.

\bibitem{AWA}
Wolfram~Research{,} Inc.
\newblock {\em Mathematica, {V}ersion 11.3}.
\newblock Wolfram Research{,} Inc., 2018.
\newblock Champaign, IL.

\bibitem{cp1}
B.~Jackson.
\newblock A zero-free interval for chromatic polynomials of graphs.
\newblock {\em Combinatorics, Probability and Computing}, 2(3):325--336, 1993.

\bibitem{JAN}
T.~Januario, S.~Urrutia, C.~C. Ribeiro, and D.~Werra.
\newblock Edge coloring: A natural model for sports scheduling.
\newblock {\em European Journal of Operational Research}, 254(1):1--8, 2016.

\bibitem{kagno}
I.~N. Kagno.
\newblock Linear graphs of degree $\leq 6$ and their groups.
\newblock {\em American Journal of Mathematics}, 68(3):505--520, 1946.

\bibitem{JK66}
J.~G. Kalbfleisch.
\newblock {\em Chromatic graphs and Ramsey's theorem}.
\newblock PhD Thesis. University of Waterloo, 1966.

\bibitem{TWG}
N.~Lord.
\newblock Graph theory as {I} have known it, by {W}. {T}. {T}utte ({B}ook
  review).
\newblock {\em The Mathematical Gazette}, 84(499):181--182, 2000.

\bibitem{mtrans79}
B.~D. McKay.
\newblock Transitive graphs with fewer than twenty vertices.
\newblock {\em Mathematics of Computation}, 33(147):1101--1121, 1979.

\bibitem{MK1984}
B.~D. McKay.
\newblock Combinatorial data.
\newblock \url{http://users.cecs.anu.edu.au/\~bdm/nauty/}, 1984.
\newblock Accessed: 2018-10-23.

\bibitem{MK19}
B.~D. McKay.
\newblock Combinatorial data.
\newblock \url{http://users.cecs.anu.edu.au/\~bdm/data/graphs.html}, 1984.
\newblock Accessed: 2018-10-23.

\bibitem{g6}
B.~D. McKay.
\newblock Description of graph6, sparse6 and digraph6 encodings.
\newblock \url{http://users.cecs.anu.edu.au/\~bdm/data/formats.txt}, 2015.
\newblock Accessed: 2018-10-23.

\bibitem{mcroyle}
B.~D. McKay and G.~F. Royle.
\newblock The transitive graphs with at most 26 vertices.
\newblock {\em Ars Combinatoria}, 30:161--176, 1990.

\bibitem{MS}
B.~D. McKay and E.~Spence.
\newblock Classification of regular two-graphs on 36 and 38 vertices.
\newblock {\em Australasian Journal of Combinatorics}, 24:293--300, 2001.

\bibitem{mcwha}
P.~McWha.
\newblock {\em Graduate course project}.
\newblock unpublished, Department of Combinatorics and Optimization, University
  of Waterloo, 1973.

\bibitem{HM}
H.~M{\'e}lot.
\newblock Facet defining inequalities among graph invariants: The system
  {G}ra{PH}edron.
\newblock {\em Discrete Applied Mathematics}, 156(10):1875--1891, 2008.

\bibitem{MM}
M.~Meringer.
\newblock Fast generation of regular graphs and construction of cages.
\newblock {\em Journal of Graph Theory}, 30(2):137--146, 1999.

\bibitem{MI2011}
S.~Mitchell, M.~O. Sullivan, and I.~Dunning.
\newblock {\em PuLP: a linear programming toolkit for python}.
\newblock The University of Auckland,
  \url{http://www.optimization-online.org/DB\_FILE/2011/09/3178.pdf}, 2011.
\newblock Accessed: 2019-02-25.

\bibitem{morris}
P.~A. Morris.
\newblock A catalogue of trees on $n$ nodes, $n < 14$, {M}athematical
  observations, research and other notes, {P}aper {N}o. 1 {S}t{A}
  (mimeographed).
\newblock {\em Publications of the Department of Mathematics, University of the
  West Indies}, 1971.

\bibitem{morriss}
P.~A. Morris.
\newblock Self-complementary graphs and digraphs.
\newblock {\em Mathematics of Computation}, 27:216--217, 1973.

\bibitem{ENS}
E.~Neufeld.
\newblock {\em Practical toroidality testing}.
\newblock MSc thesis. University of Victoria, 1993.

\bibitem{ENM}
E.~Neufeld and W.~Myrvold.
\newblock Practical toroidality testing.
\newblock In {\em Proceedings of the Eighth Annual ACM-SIAM Symposium on
  Discrete Algorithms}, SODA '97, pages 574--580, Philadelphia, PA, USA, 1997.
  Society for Industrial and Applied Mathematics.

\bibitem{TP2014}
T.~D.~P. Peixoto.
\newblock Graph-tool efficient network analysis.
\newblock \url{http://graph-tool.skewed.de}, 2014.
\newblock Accessed: 2019-02-01.

\bibitem{POT}
P.~Poto{\v{c}}nik.
\newblock A list of 4-valent 2-arc-transitive graphs and finite faithful
  amalgams of index (4,2).
\newblock {\em European Journal of Combinatorics}, 30(5):1323--1336, 2009.

\bibitem{PPG}
P.~Poto{\v{c}}nik, P.~Spiga, and G.~Verret.
\newblock Bounding the order of the vertex-stabiliser in 3-valent
  vertex-transitive and 4-valent arc-transitive graphs. arxiv preprint
  arxiv:1010.2546.
\newblock 2010.
\newblock \url{https://arxiv.org/abs/1010.2546}.

\bibitem{PP}
P.~Poto{\v{c}}nik, P.~Spiga, and G.~Verret.
\newblock A census of 4-valent half-arc-transitive graphs and arc-transitive
  digraphs of valence two. arxiv preprint arxiv:1310.6543.
\newblock 2013.
\newblock \url{https://arxiv.org/pdf/1310.6543.pdf}.

\bibitem{PSG}
P.~Poto{\v{c}}nik, P.~Spiga, and G.~Verret.
\newblock Cubic vertex-transitive graphs on up to 1280 vertices.
\newblock {\em Journal of Symbolic Computation}, 50:465--477, 2013.

\bibitem{rread}
R.~C. Read.
\newblock {\em The production of a catalogue of digraphs on 5 nodes}.
\newblock Report UWI/CCl, Computing Centre, University of the West Indies,
  1973.

\bibitem{RR81}
R.~C. Read.
\newblock A survey of graph generation techniques.
\newblock In {\em Combinatorial {M}athematics {VIII}}, pages 77--89, Editor: K.
  L. McAvaney. Berlin, Heidelberg, 1981.

\bibitem{RR91}
R.~C. Read and G.~F. Royle.
\newblock Chromatic roots of families of graphs.
\newblock {\em Graph Theory, Combinatorics, and Applications}, 2:1009--1029,
  1991.

\bibitem{RW1998}
R.~C. Read and R.~J. Wilson.
\newblock {\em An Atlas of Graphs}.
\newblock Oxford University Press, 1998.

\bibitem{reidd}
K.~B. Reid and E.~Brown.
\newblock Doubly regular tournaments are equivalent to skew {H}adamard
  matrices.
\newblock {\em Journal of Combinatorial Theory, Series A}, 12(3):332--338,
  1972.

\bibitem{RS}
N.~Robertson and P.~D. Seymour.
\newblock Graph minors. {III}. {P}lanar tree-width.
\newblock {\em Journal of Combinatorial Theory, Series B}, 36(1):49--64, 1984.

\bibitem{nr}
R.~A. Rossi and N.~K. Ahmed.
\newblock The network data repository with interactive graph analytics and
  visualization.
\newblock In {\em AAAI'15: Proceedings of the Twenty-Ninth AAAI Conference on
  Artificial Intelligence}, 2015.

\bibitem{roylethesis}
G.~F. Royle.
\newblock {\em Constructive {E}numeration of {G}raphs}.
\newblock PhD thesis, University of Western Australia, 1987.

\bibitem{GR2007}
G.~F. Royle.
\newblock Combinatorial {C}atalogues.
\newblock \url{http://staffhome.ecm.uwa.edu.au/\~00013890/}, 2007.
\newblock Accessed: 2018-10-23.

\bibitem{SC1971}
C.~C. Sims.
\newblock Computation with {P}ermutation {G}roups.
\newblock In {\em Proceedings of the second ACM symposium on Symbolic and
  algebraic manipulation}, SYMSAC'71, pages 23--28, New York,USA, 1971.

\bibitem{SIO2010}
M.~S. Sio.
\newblock {\em Towards an Online Graph Atlas for Graph Theory}.
\newblock BCompSc Hons dissertation. Clayton School of Information Technology,
  Monash University, 2010.

\bibitem{snark}
P.~G. Tait.
\newblock Remarks on the colouring of maps.
\newblock {\em Proceedings of the Royal Society of Edinburgh},
  10(729):501--503, 1880.

\bibitem{sagemath}
{The Sage Developers}.
\newblock {\em {S}ageMath, the {S}age {M}athematics {S}oftware {S}ystem
  ({V}ersion 9.0)}, 2020.
\newblock {\tt https://www.sagemath.org}.

\bibitem{shiva}
A.~D. Thomas and G.~V. Wood.
\newblock {\em Group {T}ables}.
\newblock Shiva Publishing Ltd., Orpington, 1980.

\bibitem{cp2}
C.~Thomassen.
\newblock The zero-free intervals for chromatic polynomials of graphs.
\newblock {\em Combinatorics, Probability and Computing}, 6(4):497--506, 1997.

\bibitem{cututte}
W.~T. Tutte.
\newblock A family of cubical graphs.
\newblock {\em Mathematical {P}roceedings of the {C}ambridge {P}hilosophical
  {S}ociety}, 43(4):459--474, 1947.

\bibitem{wttu}
W.~T. Tutte.
\newblock A contribution to the theory of chromatic polynomials.
\newblock {\em Canadian Journal of Mathematics}, 6(80-91):3--4, 1954.

\bibitem{tuttesy}
W.~T. Tutte.
\newblock On the symmetry of cubic graphs.
\newblock {\em Canadian Journal of Mathematics}, 11:621--624, 1959.

\bibitem{tu66}
W.~T. Tutte.
\newblock On the algebraic theory of graph colorings.
\newblock {\em Journal of combinatorial theory}, 1(1):15--50, 1966.

\bibitem{coex}
E.~W. Weisstein.
\newblock Coxeter graph.
\newblock \url{http://mathworld.wolfram.com/CoxeterGraph.html}.
\newblock Accessed: 2019-01-19.

\bibitem{Welsh}
D.~J.~A. Welsh.
\newblock {\em Complexity: Knots, Colourings and Counting}.
\newblock Cambridge University Press, USA, 1993.

\end{thebibliography}
\end{document}